\documentclass[a4paper,11pt,reqno]{article}

\usepackage{hyperref}       
\usepackage{url}            
\usepackage{amsthm}
\usepackage{amsfonts}
\usepackage{amsmath}
\usepackage{amssymb}
\usepackage{color}
\usepackage{multirow}
\usepackage{float}
\usepackage{mathtools}
\usepackage{blkarray}
\usepackage{booktabs}
\usepackage{hhline}
\usepackage{float}
\usepackage{caption}
\usepackage{faktor}
\usepackage{longtable}
\usepackage{array}
\usepackage[margin=2.7cm]{geometry}
\usepackage{enumerate}
\newcolumntype{C}[1]{>{\centering\arraybackslash}m{#1}}
\usepackage[normalem]{ulem}
\usepackage{comment}
\usepackage{fancyvrb}
\usepackage{titlesec}
\usepackage{sectsty}
\usepackage{graphicx}
\usepackage{subfigure}
\usepackage{tabularx}
\usepackage{listings}
\usepackage{xcolor}
\usepackage{colortbl} 

\usepackage{color}
\definecolor{keywordcolor}{rgb}{0.7, 0.1, 0.1}   
\definecolor{tacticcolor}{rgb}{0.0, 0.1, 0.6}    
\definecolor{commentcolor}{rgb}{0.4, 0.4, 0.4}   
\definecolor{symbolcolor}{rgb}{0.0, 0.1, 0.6}    
\definecolor{sortcolor}{rgb}{0.1, 0.5, 0.1}      
\definecolor{attributecolor}{rgb}{0.7, 0.1, 0.1} 
\definecolor{light-gray}{gray}{0.95}             

\usepackage{cite}

\theoremstyle{definition}

 \usepackage{color}
\usepackage{tikz, framed}
\usetikzlibrary{arrows,shapes,matrix,decorations.pathmorphing}
\usetikzlibrary{backgrounds}

\newcommand{\Mattia}[1]{{\color{red} Mattia: #1}}

\makeindex
\date{}
\begin{document}

\title{Teaching ``Foundations of Mathematics'' with the Lean Theorem Prover}
\author{Mattia Luciano Bottoni\thanks{Master Study at the Institute of Mathematics, University of Zurich, email: mattialuciano.bottoni@uzh.ch } \and Alberto S. Cattaneo \thanks{Institute of Mathematics, University of Zurich, email: cattaneo@math.uzh.ch}  \and Elif Sa\c{c}\i kara\thanks{IMSV, University of Bern, email: elif.sacikara@gmail.com}}

\maketitle

\abstract{This study aims to observe if the theorem prover Lean positively influences students' understanding of mathematical proving. To this end, we perform a pilot study concerning freshmen students at the University of Zurich (UZH). While doing so, we apply certain teaching methods and gather data from the volunteer students enrolled in the ``Foundations of Mathematics'' course. After eleven weeks of study covering some exercise questions implemented with Lean, we measure Lean students' performances in proving mathematical statements, compared to other students who are not engaged with Lean.
For this measurement, we interview five Lean and four Non-Lean students and we analyze the scores of all students in the final exam.
Finally, we check significance by performing a $t$-test for independent samples and the Mann-Whitney $U$-test.}\\


\section{Introduction}
Automated theorem provers and proof assistants have been around since the 1960s, when the first formal language \textit{Automath} was developed \cite{geuvers2009proof}. Since then, many theorem provers have come forth, such as \textit{Agda}, \textit{Isabelle} and \textit{Coq}. Each of them has its own strengths and weaknesses. For example, \textit{Agda} has no proof automation and \textit{Isabelle} has no dependent types in its underlying structure, meaning that predicate logic can only be done in a challenging way \cite{paulson1986natural}. While automated theorem provers are more vulnerable to mistakes, proof assistants are robust, but usually very difficult to use.

To take advantage of both systems, Lean was introduced in 2013 by Leonardo de Moura at Microsoft Research Redmond, and is one of the programming languages that acts as a bridge between proof assistant and automated theorem prover \cite{ref_Lean4}. It differs from \textit{Coq}, for instance, Lean uses a more mathematician-friendly syntax and is supported by a large community of mathematicians \cite{ref_zulipchat}, while \textit{Coq} is still used more in computer sciences circles 
\cite{ref_LoC1, ref_LoC2, ref_LoC3}. Still, both are based on adaptations of the \textit{calculus of constructions}.

Even though Lean is based on such a high level of abstract mathematics, it has gained much popularity in teaching undergraduate students in recent years, and even studies about the effects of teaching mathematics with Lean are coming up \cite{avigad2019learning, thoma2021learning}. While trends show that teaching with Lean could have positive learning effects, it remains a challenge to integrate teaching Lean into current university curricula.

This paper aims to provide further evidence of the positive effects of teaching with Lean concerning freshman students' performance in the ``Foundations of Mathematics'' course taught at the Institute of Mathematics of the University of Zurich. To this end, we learn Lean ourselves by implementing seven Lean exercise sheets based on the content of the course and we organise eleven sessions over one semester in which volunteer students are taught the foundations of mathematics with Lean.

To be able to carry out such a study, we first consider certain teaching methods like 
writing appealing goals for the students \cite{mager1978lernziele}, trying to turn extrinsic motivation into intrinsic motivation \cite{deci1993selbstbestimmungstheorie}, varying a lot between frontal teaching and experimented learning \cite{wahl2013lernumgebungen}, choosing an appropriate difficulty level \cite{vygotsky1978mind}, differentiating teaching for different skill levels \cite{wodzinski2014leistungsheterogenitat}, embracing digital teaching resources \cite{petko2014einfuhrung} and creating a learn-efficient class climate \cite{caduff2014unterrichten}.

Secondly, we consider the data gathering. More explicitly, we follow the structure of Thoma and Iannone by conducting interviews with Lean and Non-Lean students to compare their proof writing and by checking whether Lean students perform better in the final exam \cite{thoma2021learning}. To measure the significance of the results, we perform two statistical tests, the $t$-test for independent samples and the Mann-Whitney $U$-test \cite{ref_manntest}. We additionally use a questionnaire with open and closed questions to better understand the motivation and contentment of the Lean students. It is important to mention that the names of the students considered are all anonymised using Romansh language starting with \textbf{L} for Lean students and \textbf{N} for Non-Lean students. The gender of the names does not represent the gender of the actual students.

The structure of this study is as follows. Section \ref{sec:Preliminaries} provides the necessary background on the Lean interface and tactics, essential for interacting with Lean. Section \ref{sec:methods} discusses the teaching methods used in meetings with volunteer students and explains the data collection process. The results, supported by sample student solutions and statistical tests, are presented in Section \ref{sec:results}. Specifically, we evaluate students' proof structure, exam scores, and progress. Finally, Section \ref{sec:discussion} discusses our experiences in learning and teaching (with) Lean.

\section{Preliminaries}\label{sec:Preliminaries}

In this section, we briefly introduce interacting with Lean through its interface, focusing on the use of \textit{tactics}. For an extended version, and importantly, to understand the connection between tactics and \textit{type theory}, we refer the reader to \cite{mattiamaster}. Let us first present the interface of Lean as a first step to use it.

\subsection{Lean Interface}\label{sub:Leaninterface}

Several digital platforms, including online options, are available for programming with Lean. For this study, Lean was downloaded and installed for Visual Studio Code (VS Code), \cite{ref_installLean, ref_GitHub}. Once this setup is completed, a Lean project could be initiated, which, when opened in VS Code, appears as follows:
\begin{figure}[H]
\centering
    \subfigure[Lean interface at line 15.]{\includegraphics[width=0.6\textwidth]{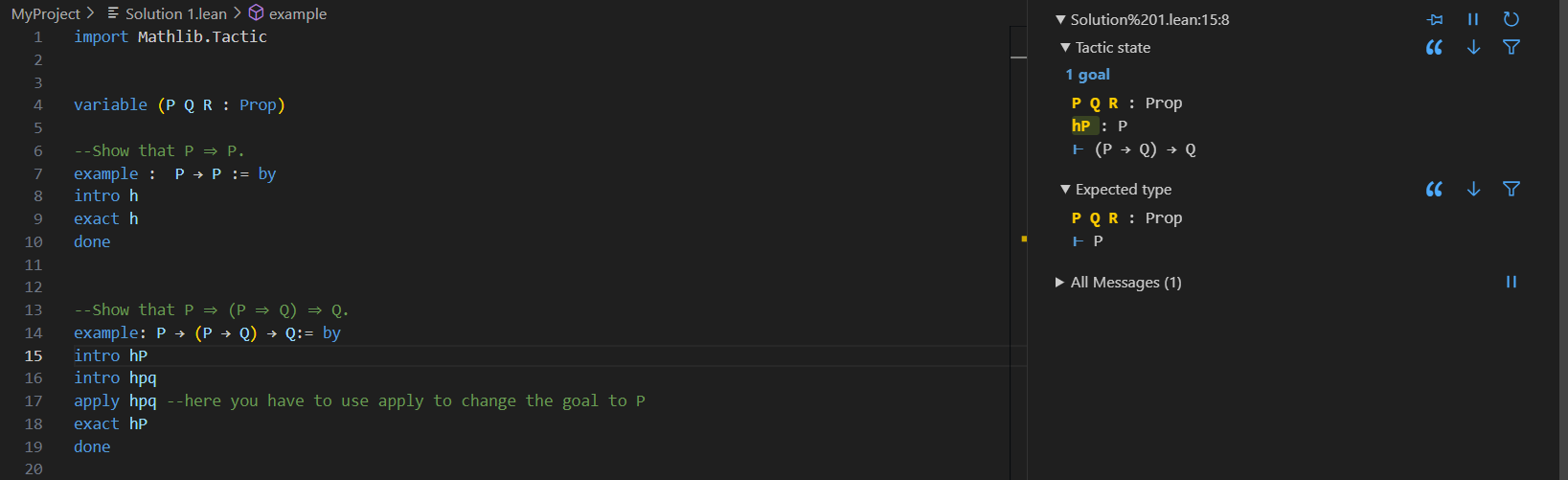}}
    \subfigure[Lean interface at line 16.]{\includegraphics[width=0.6\textwidth]{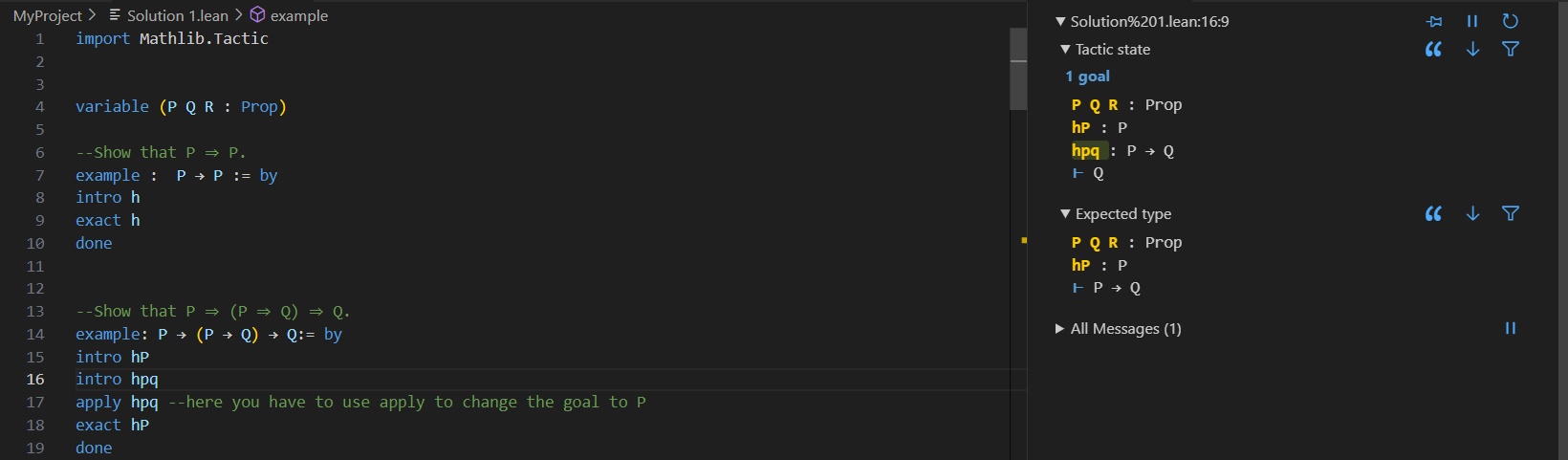}}
\caption{The Lean interface in VS Code with the helpful infoview on the right side.}
\label{fig:Lean interface}
\end{figure}

 Proof writing is done on the left side of the interface. At the top, in both figures, packages were imported to enable tactics mode and other features, which was discussed in more detail in \cite{mattiamaster}. Below, still in both figures of Figure \ref{fig:Lean interface}, an example of a short proof is shown. Note how a proposition is structured and the way the subsequent proof is constructed.

On the right side, Lean’s operations are displayed in the so-called \textit{infoview}, where the state of the proof is shown. In Figure \ref{fig:Lean interface}, depending on the line selected on the left, the infoview on the right displays a corresponding proof state. This feature is very helpful in comparison to doing proofs by hand, where hypotheses and pending conclusions must constantly be checked. If a mistake is made or a tactic not allowed is attempted, an error message appears in the infoview, as illustrated in Figure \ref{fig:error message} below.

\begin{figure}[H]
\centering
\includegraphics[width=0.4\textwidth]{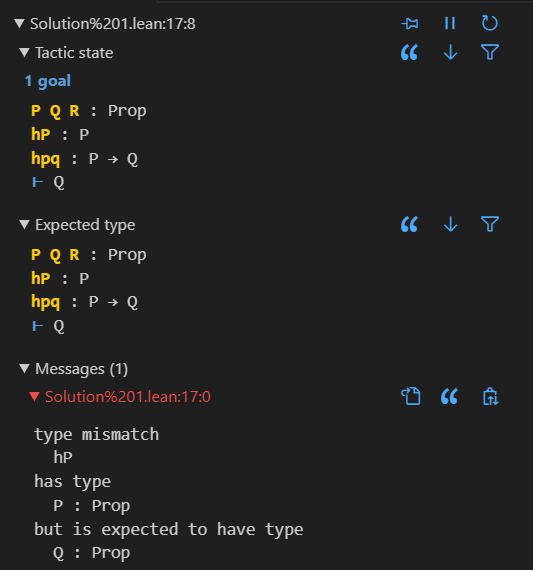}
\caption{Error-message when trying to exact $hP$ from the example above.}
\label{fig:error message}
\end{figure}

 The error message not only indicates that a mistake has been made but also specifies the type of mistake, which, depending on one’s experience, can help with direct correction. Imagine if a piece of paper could offer such feedback! In certain situations, the infoview can also display proposed theorems to aid in solving a proof, a feature discussed further when Lean’s tactics are examined.

Lean’s automation helps the proof process by providing progress feedback, identifying mistakes, and even suggesting relevant theorems from Mathlib’s extensive collection \cite{ref_mathlib}. Additional assistance is available when using Lean in VS Code. With the development of AI tools like ChatGPT, a code-suggestion AI called GitHub Copilot has been created and can be used directly in VS Code. This package requires installation and login, which is available for free with a UZH account.

\begin{figure}[H]
\centering
\includegraphics[width=0.6\textwidth]{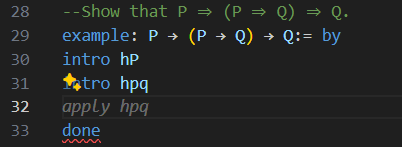}
\caption{A suggestion for the next step generated by GitHub Copilot.}
\label{fig:copilot}
\end{figure}
During the meetings with students, the copilot was intentionally disabled, as it simplified some exercises too much. However, it was still very useful for developing solutions to the exercise sheets. It should be noted that GitHub Copilot does not always provide helpful code suggestions and occasionally proposes lines that result in error messages. Nevertheless, it can assist when one is stuck at a certain point in a proof, even if only by sparking an idea through a line of code that may not work initially.

\subsection{Tactics}
When proving statements in Lean, one can use so-called \textit{tactics}, which are similar to instructions for proofs \cite{ref_Lean4}. In this section, only the tactics that positively impacted Lean students, as discussed in Section  \ref{sec:results}, are presented. In the examples shown here, we often just show a small part of the whole proof; the complete proofs can be found on GitHub \cite{ref_GitHub}. 

All tactics taught to participants throughout the semester, along with their logical foundations in \textit{natural deduction rules} and by \textit{the PAT interpretation}, can be found in \cite{mattiamaster} and the references therein. More specifically, some of the tactics presented there are directly related to derivation rules from the \textit{calculus of constructions}, used to build a term (the proof) corresponding to the given type (the proposition). For these tactics, we always present the natural deduction rule the tactic corresponds to, the proof using the tactics, and the proof in \textit{term-style}. Term-style proofs are another way to prove statements in Lean, and they are closely related to the calculus of constructions. One can even use both methods combined to build a proof. Apart from tactics and derivation rules, Lean also has access to countless theorems and definitions implemented in a huge library called Mathlib \cite{ref_mathinLean4}. We have summarized the most important addition, subtraction, multiplication, and division theorems in a ``cheat sheet" that can be found on \cite{ref_GitHub}. We note that, although we prepared such documentation in \cite{mattiamaster}, we did not consider this high-level logic in our studies with students. Therefore, we prefer not to present it here.

\subsubsection*{intro tactic}\label{sub:intro}

One of the first things to do in a proof is to introduce the variables and hypotheses that are available. This can be done when there is an implication type in the goal In Lean, the \textbf{intro} or \textbf{intros} tactics are used for introducing these variables or hypotheses.

What we are doing here is introducing $P$ and from $P$ we want to derive $Q$. This then gives the proof that $P$ implies $Q$, seen below the line.

    \begin{figure}[H]
\centering
    \subfigure[\textit{intro} used in VS Code.]{\includegraphics[width=0.45\textwidth]{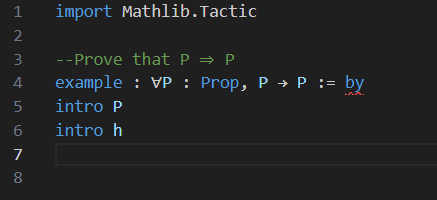}}
    \hspace{0.05\textwidth}
    \subfigure[Lean infoview after using \textit{intro}.]{\includegraphics[width=0.45\textwidth]{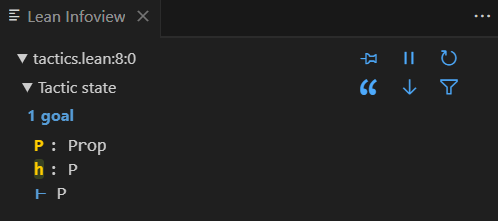}}
\caption{Using the \textit{intro} tactic to prove that $P \Rightarrow P$ holds for any proposition.}
\label{fig:intro tactic}
\end{figure}

We see in the Figure \ref{fig:intro tactic} above, on the left, that the original goal is to prove that $\forall P : \mathtt{Prop}, ~ P \rightarrow P$. Here, $P : \mathtt{Prop}$ means that $P$ is of type $\mathtt{Prop}$ (propositions). In the right picture, we see the remaining goal ($\vdash P$) and the two hypotheses we now have. They are $P : \mathtt{Prop}$, i.e.\ we now have an arbitrary proposition and $h : P$. The second means that we have a proof $h$ for $P$ independent of what proposition $P$ is exactly.

\subsubsection*{exact tactic}\label{sub:exact}

To complete the proof in Figure \ref{fig:intro tactic}, the hypothesis $h : P$ needs to be used. This can be accomplished with the \textbf{exact} tactic, which applies a proposition or proof to satisfy the given goal.

    \begin{figure}[H]
\centering
    \subfigure[\textit{exact} used in VS Code.]{\includegraphics[width=0.45\textwidth]{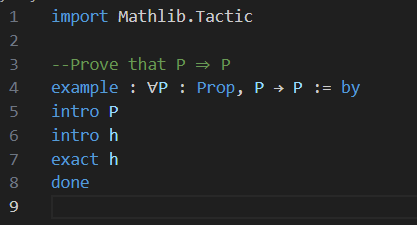}}
    \hspace{0.05\textwidth}
    \subfigure[Lean infoview after using \textit{exact}.]{\includegraphics[width=0.45\textwidth]{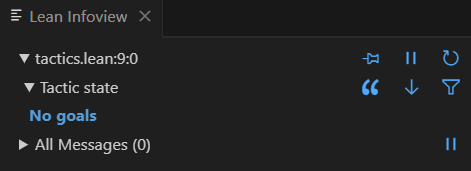}}
\caption{Using the \textit{exact} tactic concludes the proof that $P \Rightarrow P$ holds for any proposition.}
\label{fig:exact tactic}
\end{figure}

When we try to prove this without tactics, we do it as follows.
    \begin{lstlisting}[backgroundcolor = \color{light-gray}]
example : P → P := fun p : P => p
\end{lstlisting}

Here we introduce a function \textbf{fun} from $P$ to $P$ called $p$. When we are not in Lean's tactic mode, this already finishes the proof.

\subsubsection*{apply tactic}\label{sub:apply}

Unfortunately it is not always straightforward to prove a statement. For example, Consider a case with two implications in the statement to prove, such as $P \Rightarrow (P \Rightarrow Q) \Rightarrow Q$. In this situation, simply using the \textit{intro} tactic to introduce $P$ and $P \Rightarrow Q$ and then applying \textit{exact} will not work, as \textit{exact} requires the hypothesis and goal to be identical. When the hypothesis contains an implication and the goal aligns with the conclusion of that implication, the \textbf{apply} tactic is suitable.

\begin{figure}[H]
\centering
    \subfigure[Before we use \textit{apply}.]{\includegraphics[width=0.45\textwidth]{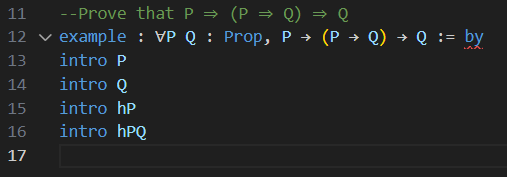}}
    \hspace{0.05\textwidth}
    \subfigure[Lean infoview before using \textit{apply}.]{\includegraphics[width=0.45\textwidth]{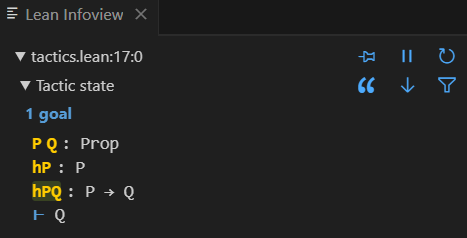}}
    \subfigure[\textit{apply} used in VS Code.]{\includegraphics[width=0.45\textwidth]{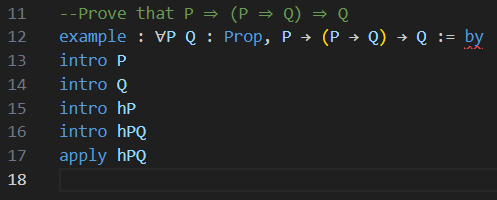}}
    \hspace{0.05\textwidth}
    \subfigure[Lean infoview after using \textit{apply}.]{\includegraphics[width=0.45\textwidth]{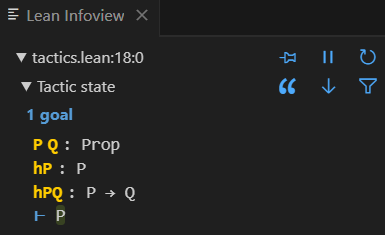}}
\caption{Using the \textit{apply} tactic to deal with hypotheses including implications.}
\label{fig:apply tactic}
\end{figure}

The right-hand side of the hypothesis $hPQ$ is $Q$, which is the same as the goal. This allows us to use the \textit{apply} tactic to shift the goal to the left-hand side of the hypothesis, $P$. The \textit{apply} tactic operates on this principle: given an implication, such as $P \Rightarrow Q$ as a hypothesis, with $hPQ$, as the proof of the implication, it suffices to show the left-hand side $P$ to validate the right-hand side $Q$. Thus, we can change the goal from the right-hand side to the left-hand side. We note that, alternatively, one can forget the \textit{apply} tactic and use \textit{exact hPQ hP} instead.

Proving this statement without tactics would look like this.
\begin{lstlisting}[backgroundcolor = \color{light-gray}]
example : P → (P → Q) → Q :=
  fun p : P =>
  fun f : P → Q =>
  f p
\end{lstlisting}

One may have already noticed that the structure of the proof changes when we construct it without tactics. Tactics can only be used inside a by-done environment. In this case, we see that \textit{apply} was not used explicitly.

\subsubsection*{rw tactic}\label{sub:rw}

We are somewhat more fortunate when we have if-and-only-if statements as hypotheses, as this allows us to decide which side to use. When an if-and-only-if statement appears as a hypothesis, or in Mathlib, we can apply the \textbf{rw} (rewrite) tactic to adjust the statement accordingly. The \textbf{rw at} tactic can even be used to modify a specific hypothesis. If multiple rewrites are required, we can use \textbf{repeat rw}.

    \begin{figure}[H]\ContinuedFloat
\centering
    \subfigure[Before we use \textit{rw}.]{\includegraphics[width=0.35\textwidth]{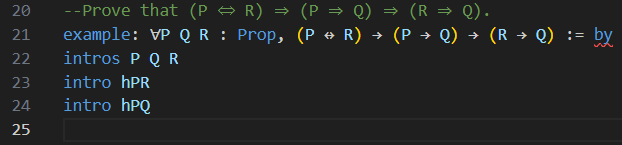}}
    \hspace{0.05\textwidth}
    \subfigure[Lean infoview before using \textit{rw}.]{\includegraphics[width=0.35\textwidth]{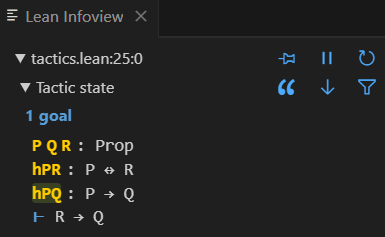}}
    \subfigure[\textit{rw} used in the goal.]{\includegraphics[width=0.35\textwidth]{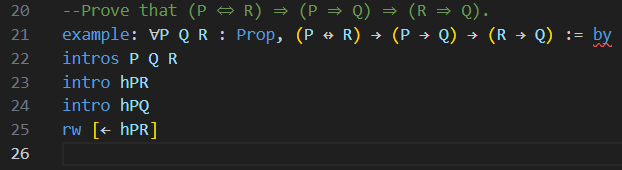}}
    \hspace{0.05\textwidth}
    \subfigure[Lean infoview after using \textit{rw} in the goal.]{\includegraphics[width=0.35\textwidth]{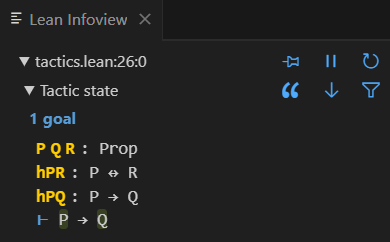}}
    \end{figure}
    \begin{figure}[H]
    \centering
    \subfigure[\textit{rw} used in the hypothesis.]{\includegraphics[width=0.35\textwidth]{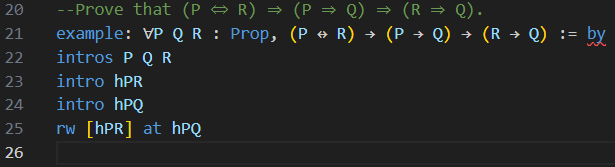}}
    \hspace{0.05\textwidth}
    \subfigure[Lean infoview after using \textit{rw} in the hypothesis.]{\includegraphics[width=0.35\textwidth]{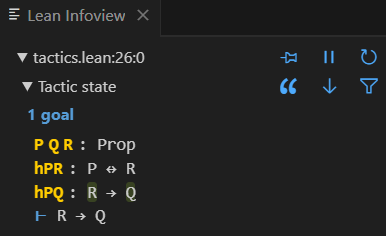}}
\caption{Using the \textit{rw} tactic to change between $P$ and $R$.}
\label{fig:rw tactic}
\end{figure}

In Figure \ref{fig:rw tactic}, we illustrate the difference between using the \textit{rw} tactic in the goal versus in a hypothesis. For the first case, a left arrow ($\leftarrow$) is required, as we want to rewrite $R$ to $P$, not vice versa. In this specific example, changing the hypothesis is not ideal, as the label $hPQ$ may lose its meaning after $P$ is rewritten as $R$. Alternatively, we could assign a more general label to the hypothesis, avoiding direct references to $P$ or $R$.

\subsubsection*{constructor tactic}\label{sub:constructor}

If-and-only-if statements frequently appear that we would like to prove. On paper, this is typically done by proving each implication separately. In Lean, the same approach applies. The \textbf{constructor} tactic is used here to split an if-and-only-if goal into two subgoals, requiring both implications to be proved.

\begin{figure}[H]\ContinuedFloat
\centering
    \subfigure[Before we use \textit{constructor}.]{\includegraphics[width=0.35\textwidth]{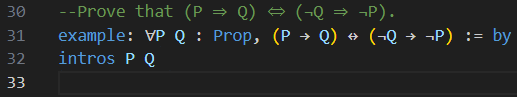}}
    \hspace{0.05\textwidth}
    \subfigure[Lean infoview before using \textit{constructor}.]{\includegraphics[width=0.35\textwidth]{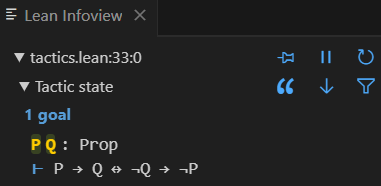}}
    \end{figure}
    \begin{figure}[H]
    \centering
    \subfigure[\textit{constructor} used in VS Code.]{\includegraphics[width=0.35\textwidth]{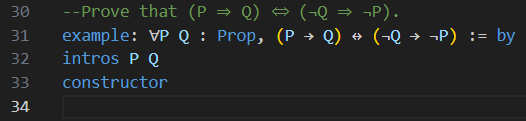}}
    \hspace{0.05\textwidth}
    \subfigure[Lean infoview after using \textit{constructor}.]{\includegraphics[width=0.35\textwidth]{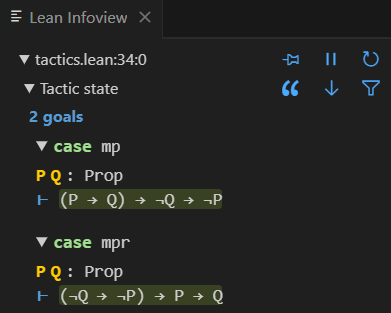}}
\caption{The \textit{constructor} tactic splits a goal into two subgoals for both implications.}
\label{fig:constructor tactic}
\end{figure}    

The \textit{constructor} tactic can also be used when we have a $\wedge$ in the goal.

    \begin{figure}[H]\ContinuedFloat
\centering
    \subfigure[Before we use \textit{constructor}.]{\includegraphics[width=0.45\textwidth]{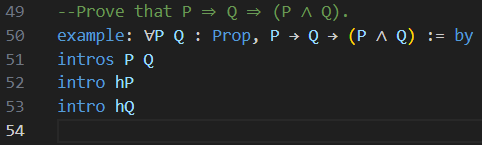}}
    \hspace{0.05\textwidth}
    \subfigure[Lean infoview before using \textit{constructor}.]{\includegraphics[width=0.45\textwidth]{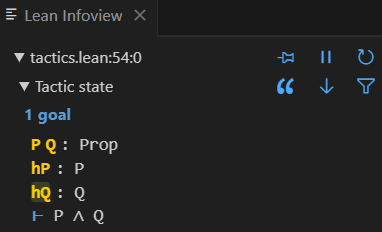}}
    \end{figure}
    \begin{figure}[H]
    \centering
    \subfigure[\textit{constructor} used in VS Code.]{\includegraphics[width=0.45\textwidth]{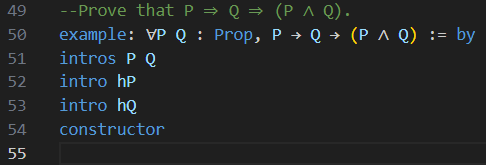}}
    \hspace{0.05\textwidth}
    \subfigure[Lean infoview after using \textit{constructor}.]{\includegraphics[width=0.45\textwidth]{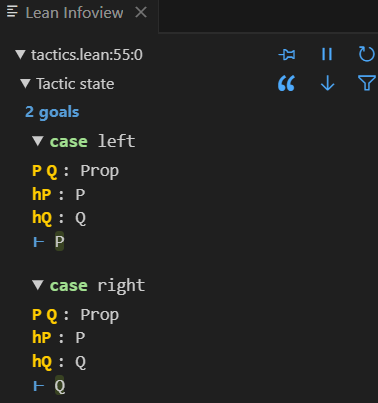}}
\caption{The \textit{constructor} tactic splits a goal into two subgoals for both sides of the $\wedge$ operator.}
\label{fig:and constructor tactic}
\end{figure}

After we split the proof into two subgoals, we can simply use exact $hP$ and exact $hQ$ to finish the proof. When we are not in the tactics mode, we use the same expression as in the sequent calculus:
\begin{lstlisting}[backgroundcolor = \color{light-gray}]
example : P → Q → (P ∧ Q) :=
  fun p =>
  fun q =>
  And.intro p q
\end{lstlisting}

With the \textit{constructor} tactic, a major advantage of tactics becomes clear. While the $\wedge$-introduction rule works only for conjunction terms, the \textit{constructor} tactic can also split if-and-only-if statements. Lean automatically recognizes which case is being handled.

\subsubsection*{cases' tactic}\label{sub:cases}

Disjunction ($\vee$) and conjunction ($\wedge$) terms can also appear as hypotheses. Consider how they differ when the statements are split. Whenever there is a $\wedge$ or a $\vee$ in one of the hypotheses instead of in the goal, we can use the \textbf{cases'} tactic. For a logical ``and'', this produces a single proof term with two hypotheses, whereas for a logical ``or'', the statement must be proven twice: once with the left-hand side of the hypothesis and once with the right-hand side.

    \begin{figure}[H]
\centering
    \subfigure[Before we use \textit{cases'}.]{\includegraphics[width=0.35\textwidth]{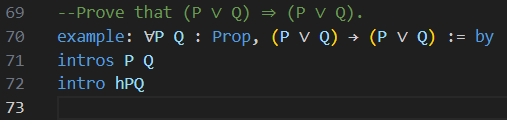}}
    \hspace{0.05\textwidth}
    \subfigure[Lean infoview before using \textit{cases'}.]{\includegraphics[width=0.35\textwidth]{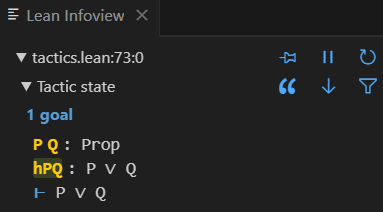}}
    \subfigure[\textit{cases'} used in VS Code.]{\includegraphics[width=0.35\textwidth]{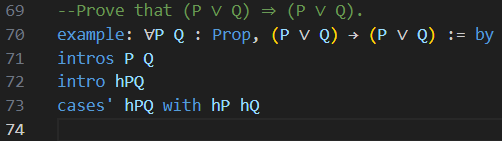}}
    \hspace{0.05\textwidth}
    \subfigure[Lean infoview after using \textit{cases'}.]{\includegraphics[width=0.35\textwidth]{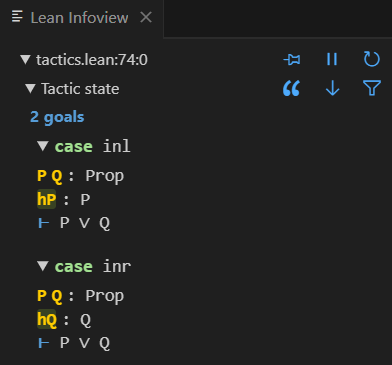}}
\caption{\textit{cases'} applied on a hypothesis with ``or'' gives us two subgoals to prove.}
\label{fig:cases' or tactic}
\end{figure}

Using the $\wedge$-elimination rule, we can derive both $P$ and $Q$. When applying the \textit{cases'} tactic on a conjunction hypothesis, Lean performs both derivations, resulting in two separate hypotheses. The elimination rule for disjunction is somewhat more complex. Assume that we have a statement $R$ and a proof of $P \vee Q$. To show that $P \vee Q$ implies $R$, we need to prove both that $P$ implies $R$ and that $Q$ implies $R$. Each must imply $R$, as we want the flexibility to use either $P$ or $Q$ to derive $R$. Consequently, when we encounter $P \vee Q$ in a hypothesis, we must establish two statements: $P \Rightarrow R$ and $Q \Rightarrow R$.  

Lean can determine from the context which of the two deduction rules should be applied. Of course, we could use \textit{exact} $hPQ$ directly; however, this example demonstrates the impact of the \textit{cases'} tactic on a logical or-statement. Note that we can assign individual names to the two newly created hypotheses.

Without using the tactics mode, it is essential to distinguish between conjunctions and disjunctions in the proof. For disjunction, the proof proceeds as follows:
\begin{lstlisting}[backgroundcolor = \color{light-gray}]
example : (P ∨ Q) → (Q ∨ P) :=
  fun h : P ∨ Q =>
  Or.elim h
    (fun p : P => Or.inr p)
    (fun q : Q => Or.inl q)
\end{lstlisting}

A significant advantage of the tactic mode is to show current goal. This feature is especially helpful in cases where two independent proofs are required, as illustrated in the example above.

    \begin{figure}[H]
\centering
    \subfigure[Before we use \textit{cases'}.]{\includegraphics[width=0.45\textwidth]{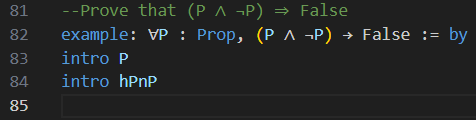}}
    \hspace{0.05\textwidth}
    \subfigure[Lean infoview before using \textit{cases'}.]{\includegraphics[width=0.45\textwidth]{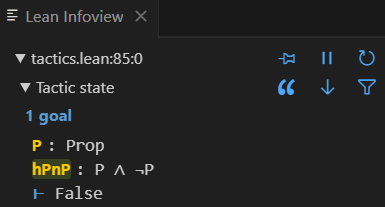}}
    \subfigure[\textit{cases'} used in VS Code.]{\includegraphics[width=0.45\textwidth]{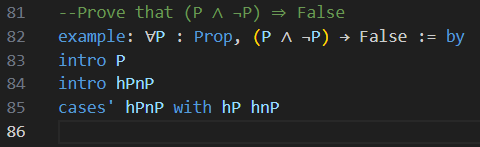}}
    \hspace{0.05\textwidth}
    \subfigure[Lean infoview after using \textit{cases'}.]{\includegraphics[width=0.45\textwidth]{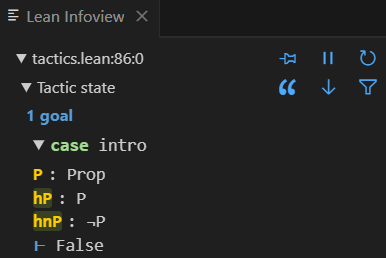}}
\caption{\textit{cases'} applied on a hypothesis with ``and'' gives us two hypotheses to prove one goal.}
\label{fig:cases' and tactic}
\end{figure}

Notice that, in this case, the goal remains unchanged when using the cases' tactic. Keep in mind that $\neg P$ and $P \Rightarrow \text{False}$ are logically equivalent; this can be verified using a truth table for the bachelor level students. In Lean, $\neg P$ is defined as $P \Rightarrow \text{False}$ and can be used accordingly. To complete the proof, we \textit{apply} $hnP$, which changes the goal to $P$, and then use \textit{exact} $hP$.

We can do the same without tactics. Notice that here we do not use \textit{apply} again:
\begin{lstlisting}[backgroundcolor = \color{light-gray}]
example : (P ∧ ¬P) → False :=
  fun h : P ∧ ¬P =>
  And.right h (And.left h)
\end{lstlisting}
While it is convenient that the rules are named \textit{And.right} and \textit{And.left}, it should be noted that, unlike when using tactics, we do not see the two hypotheses $hP$ and $hnP$ in this case.

\subsubsection*{symm tactic}\label{sub:symm}

Lean is very precise when checking proofs. For example, if the hypothesis is $h : x = 3$, we cannot prove the goal $\vdash 3 = x$ using exact $h$. One needs to switch the hypothesis or the goal. The \textbf{symm} tactic lets us interchange the left- and right-hand side of an equality in any goal or hypothesis.

    \begin{figure}[H]\ContinuedFloat
\centering
    \subfigure[Before we use \textit{symm}.]{\includegraphics[width=0.45\textwidth]{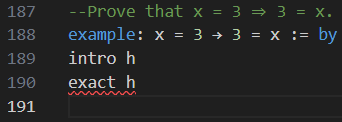}}
    \hspace{0.05\textwidth}
    \subfigure[Lean infoview before \textit{symm}.]{\includegraphics[width=0.45\textwidth]{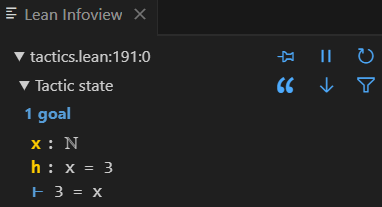}}
    \end{figure}
    \begin{figure}[H]
    \centering
    \subfigure[\textit{symm} used in VS Code.]{\includegraphics[width=0.45\textwidth]{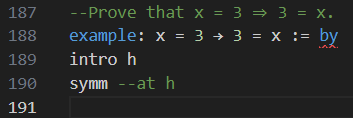}}
    \hspace{0.05\textwidth}
    \subfigure[Lean infoview after using \textit{symm}.]{\includegraphics[width=0.45\textwidth]{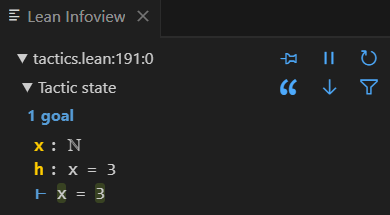}}
\caption{Without \textit{symm}, the proof cannot be concluded.}
\label{fig:symm tactic}
\end{figure}

We could also change the equality in the hypothesis $h$, then the new hypothesis would be $h : 3 = x$. Exactly like the goal.

\subsubsection*{induction' tactic}\label{sub:inductiontactic}

Natural induction is built into Lean’s capabilities. Using the \textbf{induction'} tactic, proofs by induction can be initiated over a variable, with options available to perform proofs by strong induction as well.

    \begin{figure}[H]
\centering
    \subfigure[A proof using natural induction.]{\includegraphics[width=0.45\textwidth]{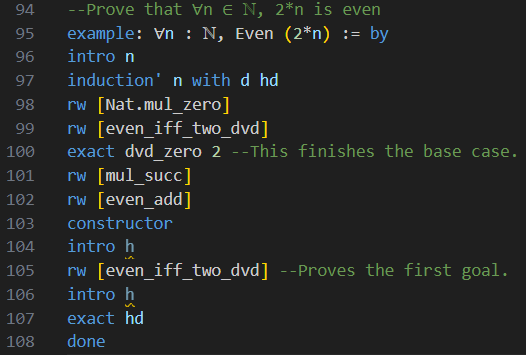}}
    \hspace{0.05\textwidth}
    \subfigure[Lean infoview after using \textit{induction'}.]{\includegraphics[width=0.45\textwidth]{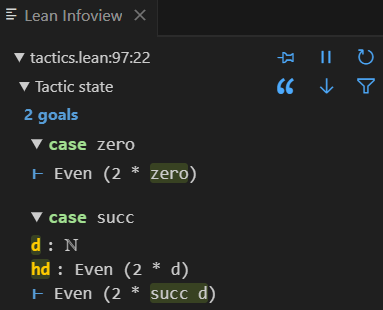}}
\caption{Natural induction in Lean.}
\label{fig:induction tactic}
\end{figure}

\noindent Even though this is a trivial statement, we use quite specific theorems to prove it. We elaborate on that in Chapter \ref{subsec:proconlearning}.\\

\noindent The proof by the strong induction method was used to solve a question from one of the exercise sheets. Interested readers can have a look at that proof in Level 3, Exercise 3.2 in \cite{ref_GitHub}.\\

\subsubsection*{ring\_nf tactic}\label{sub:ringnf}
  As long as our goal contains expressions with only ring arithmetic, we can use the \textbf{ring\_nf} tactic to prove the statement.
    \begin{lstlisting}[backgroundcolor = \color{light-gray}]
--Prove that (x+y)^2 = x^2 + 2xy + y^2.
example: ∀x y : ℝ, (x+y)^2 = x^2 + 2*x*y + y^2 := by
intros x y
ring_nf
done
\end{lstlisting}

\section{Methods}\label{sec:methods}
In this section, we consider the teaching methods and data gathering process. The teaching methods were used to prepare the weekly sessions we held with five volunteer first-year mathematics students.
\subsection{Teaching Methods}\label{subsec:teachingmethods}

\begin{figure}[H]
\centering
\includegraphics[width=0.8\textwidth]{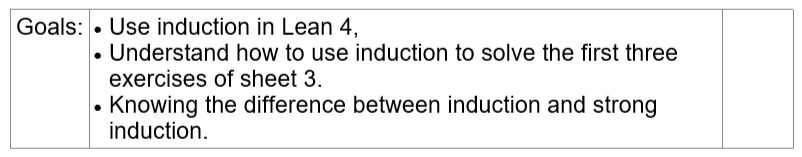}
\caption{The goals from our 4-th meeting.}
\label{fig:goalsexample}
\end{figure}

At the start of each meeting, we showed students the session goals to help them focus on key topics. This approach has been shown to improve student achievement \cite{frey2010ausgewahlte}. The goals were structured according to Mager's method \cite{mager1978lernziele}, including both content and the skills students should develop. Figure \ref{fig:goalsexample} provides an example from one of the meetings.

  We distinguish between extrinsic and intrinsic motivation. Extrinsic motivation comes from societal goals and expectations, while intrinsic motivation arises from personal interests and curiosity. Deci and Ryan argued in \cite{deci1993selbstbestimmungstheorie} that these two types are not opposites; extrinsic motivation can become intrinsic through a four-stage process: external regulation, introjected regulation, identified regulation, and integrated regulation. This process describes how external motivation becomes internalized, leading to self-determination and intrinsic motivation. With the questionnaire in Section \ref{subsec:interviews}, we aim to assess the motivation of Lean students to determine if higher motivation leads to better learning outcomes.

The meetings followed the sandwich method \cite{wahl2013lernumgebungen}, which alternates between direct instruction and hands-on activities and allows students to work independently for a significant portion of the lesson. On the other hand, when selecting exercises to discuss, we considered an appropriate difficulty level. If the content is too difficult, students become overwhelmed and frustrated; if it is too easy, they get bored \cite{Petko22}. We aimed to reach the `zone of proximal development' \cite{vygotsky1978mind}, where the difficulty is just beyond students' current knowledge but still attainable. This concept influenced the restructuring of the meetings in the final weeks of the semester, as explained in Section \ref{subsec:studygroups}.

Heterogeneity is a common topic in teaching \cite{wodzinski2014leistungsheterogenitat}, as classrooms often include a wide range of skill levels. While individual attention for each student is ideal, it is often unfeasible in large classes. However, with our small group, we were able to tailor instruction to each student's needs. For instance, one of the Lean student, Ladina, had more advanced knowledge compared to the others, and even us. The extra one-on-one meetings with her enabled us to provide differentiated teaching \cite{hofer2021individuelle}, allowing us to increase the difficulty level specifically for her.

  The need to adapt our teaching to the rise of digital media is clear, yet change has been slow, and effective tools for integrating digital technologies into classrooms are still lacking \cite{petko2014einfuhrung}. Although digital technology has been present for nearly 40 years, research on its effective use in education takes time. Some methods are emerging, emphasizing the importance of digital teaching. Schools must adapt as the world evolves rapidly. Research highlights the need for media-aware teaching, which fosters critical thinking and the ability to manage information flow \cite{hughes2016guiding}, while others caution that digital media should enhance, not replace, good teaching \cite{freeman2017nmc}.

In our case, we combined both methods mentioned above, teaching digital skills by integrating digital resources without losing the core content. We are not alone in this approach. The term STEM (Science, Technology, Engineering, and Mathematics) has gained prominence in recent years. Both \cite{hughes2016guiding} and \cite{freeman2017nmc} emphasize the importance of integrating multiple subjects, and Lean effectively combines technology and mathematics. Additionally, Lean enables students to develop information literacy, critical thinking, creative digital work, and interdisciplinary STEM skills, which is not something to be taken for granted.

  We maintained an informal, relationship of equals during sessions, which helped build a strong connection with the students and encouraged their participation in interviews. However, it's important that a good atmosphere does not compromise learning, as students are still learners, not peers \cite{caduff2014unterrichten}. To maintain this balance, we arranged the classroom so that students sat in a half-circle, fostering closeness without reducing formality. A supportive class climate not only enhances learning but also boosts self-confidence, motivation, positive attitudes toward lessons and school, as well as social behavior and interest development \cite{meyer2014guter}. This, in turn, promotes intrinsic motivation.
  
\subsection{Data Gatheting}\label{subsec:datagathering}

The main goal of this study is to observe the effect of both teaching Lean and teaching with Lean on students' proving skills. To achieve this, we need to find a way to measure the quality of proofs and design appropriate interviews and questionnaires for both Lean and Non-Lean students.

\subsection*{Proof Structure}\label{subsec:proofstructure}
Thoma and Iannone outline various criteria for assessing the quality of a proof \cite{thoma2021learning}, and we use these to evaluate the students' performance. The grading tables are available in the Interviews and Questionnaire section, with the results presented in Section \ref{subsec:interviewresults}.

\begin{itemize}
    \item \textbf{Definitions and their use}: Check whether students understand the given definitions and can correctly apply them in their proofs.
    \item \textbf{Mathematical symbols and their use}: Check whether students correctly understand and apply mathematical symbols, which are essential for constructing a valid proof.
    \item \textbf{Logical statements and their links}: Check whether students understand how to deal with implications or equivalences, and whether they correctly handle quantifiers.
    \item \textbf{High level ideas}: Check whether students grasp the key idea required for a proof and can find links between different mathematical theories, such as set theory and first-order logic.
    \item \textbf{Modular structure of the proof}:  Check the ability of students to break a proof into different sub-proofs and to combine them correctly at the end, ensuring the overall structure is solid.
    \item \textbf{Use of examples}: Check whether students can use examples effectively to find patterns that might help them complete the proof.
\end{itemize}

\subsection*{Interviews and Questionnaire} \label{subsec:interviews}
We interviewed both groups, Lean and Non-Lean students, through a quantitative test. To gain further insights into the time spent with Lean, we also prepared a separate questionnaire exclusively for Lean students. In this section, we present and justify the questions selected for the interview, discuss the main focuses of each question, and outline the grading tables used to assess the interviews based on the grading criteria explained above. Finally, we present the questionnaire. The results of the interviews and the evaluation of the questionnaire can be found in Section \ref{sec:results}.

The four interview questions are the following.
\begin{enumerate}
    \item \textit{Prove the following deMorgan law without a truth table:}
    $$\neg(p \vee q) = (\neg p)\wedge (\neg q)$$ 
    With this question we aim to check if Lean students developed additional skills to deal with propositional logic proofs and the way that Non-Lean students approach these kinds of proofs.
    \item \textit{Using all axioms of addition and multiplication of natural numbers, prove:}
    $$\forall a, b, c \in \mathbb{N}, ~ a \cdot (b + c) = a \cdot b + a \cdot c $$ In this question, we want to measure the understanding of proof by induction and the definition of the successor function. We also evaluate their ability to handle propositional equalities correctly in the base case.
    \item \textit{Recall that:} 
    $$a \equiv b \mod n \Leftrightarrow (\exists c, d, r \in \mathbb{N}, ~ a = c \cdot n + r  ~ \mbox{and} ~ b = d \cdot n + r) \Leftrightarrow (\exists k \in \mathbb{Z}, a - b = k \cdot n).$$
    \textit{Show that} congruence mod $n$ \textit{with} $n = 3$ \textit{is an equivalence relation.}\\ 
    With this question we aim to observe if students can apply the proof of an equivalence relation in the context of modulo. We also evaluate their understanding of the definition of modulo. 
    \item \textit{If} $J \neq \emptyset$ \textit{and} $J \subseteq I$\textit{, does it follow that} $\bigcup_{\alpha \in J} A_{\alpha} \subseteq \bigcup_{\alpha \in I} A_{\alpha}$\textit{? What about} $\bigcap_{\alpha \in J} A_{\alpha} \subseteq \bigcap_{\alpha \in I} A_{\alpha}$\textit{?}\\~
    \\
    For this last question, we want to know more about students' understanding of set-theoretic proofs and the definitions of some more complex set structures.
\end{enumerate}~

   The idea of the interview is to cover as many topics from the lecture as possible. We see that there is a question about propositional logic, one about natural numbers and natural induction, one about the definition of an equivalence relation, and one on sets, specifically on indexed sets and operations on sets.  Table \ref{tab:overviewtopics} demonstrates how often a topic is covered in either the Thoma and Iannone paper \cite{thoma2021learning}, the interviews we conducted, the exercise sheets of the course ``Foundations of Mathematics'', and the final exam of that course. One can see that there is a high overlap between the topics chosen.

\begin{table}[H]
\centering
\resizebox{\textwidth}{!}{
\begin{tabular}{|l|l|l|l|l|l|l|l|l|}
\hline
\textbf{}                & \textbf{\begin{tabular}[c]{@{}l@{}}Propositional \\ Logic\end{tabular}} & \textbf{\begin{tabular}[c]{@{}l@{}}Basic \\ Set Theory\end{tabular}} & \textbf{Proofs} & \textbf{\begin{tabular}[c]{@{}l@{}}Mathematical \\ Induction\end{tabular}} & \textbf{\begin{tabular}[c]{@{}l@{}}Relations and \\ Functions\end{tabular}} & \textbf{Cardinality} & \textbf{\begin{tabular}[c]{@{}l@{}}Natural \\ Numbers\end{tabular}} & \textbf{\begin{tabular}[c]{@{}l@{}}Order \\ Relations\end{tabular}} \\ \hline
\textbf{\begin{tabular}[c]{@{}l@{}}Thoma and \\ Iannone\end{tabular}}           &                                                                         &                                                                      & 2            & 1                                                                          & 2                                                                        &                      & 1                                                                   &                                                                     \\ \hline
\textbf{Interview}       & 1                                                                       & 1                                                                    &                 & 1                                                                          & 1                                                                           &                      & 1                                                                   &                                                                     \\ \hline
\textbf{\begin{tabular}[c]{@{}l@{}}Exercise \\ Sheets\end{tabular}} & 2                                                                    & 2                                                                 & 1               & 3                                                                    & 2                                                                        & 1                    & 2                                                                & 1                                                                   \\ \hline
\textbf{Exam}            & 1                                                                       &                                                                      & 2            & 1                                                                          & 2                                                                        & 1                    &                                                                     & 1                                                                   \\ \hline
\end{tabular}}
\caption{Overview of the topics asked about in various examinations.}
\label{tab:overviewtopics}
\end{table}

  Question 1 is specifically designed to give the Lean students an advantage, while Question 4 is better suited for students who did not learn Lean. This is because Non-Lean students have not been taught how to prove propositional logic without truth tables, and we did not spend much time on exercises like Question 4 with the Lean students. In these four questions, in addition to the focuses mentioned above, we evaluate their proof writing based on the characteristics explained in Section \ref{subsec:proofstructure}. The following tables show the grading criteria used to assess the students' exercises.

\begin{table}[H]
\parbox{.5\linewidth}{
\centering
\resizebox{.5\textwidth}{!}{
\begin{tabular}{|llll|}
\hline
\multicolumn{4}{|l|}{\cellcolor[HTML]{E59EDC}\textbf{Question 1}}                                         \\ \hline
\multicolumn{1}{|l|}{} & \multicolumn{1}{l|}{0 Points} & \multicolumn{1}{l|}{1 Point} & 2 Points \\ \hline
\multicolumn{1}{|l|}{\textbf{D}} &
  \multicolumn{1}{l|}{\begin{tabular}[c]{@{}l@{}}Needs help with the\\ definitions, even when\\ given.\end{tabular}} &
  \multicolumn{1}{l|}{\begin{tabular}[c]{@{}l@{}}Understands the\\ needed definitions, i.e.\\ ``or''/``and'', but maybe\\ does not know them all.\end{tabular}} &
  \begin{tabular}[c]{@{}l@{}}Uses them in a\\ confident way.\end{tabular} \\ \hline
\multicolumn{1}{|l|}{\textbf{MS}} &
  \multicolumn{1}{l|}{\begin{tabular}[c]{@{}l@{}}Struggles with writing\\ mathematics.\end{tabular}} &
  \multicolumn{1}{l|}{\begin{tabular}[c]{@{}l@{}}Has a cLean way of\\ writing mathematics.\end{tabular}} &
  \begin{tabular}[c]{@{}l@{}}No neglection\\ detected.\end{tabular} \\ \hline
\multicolumn{1}{|l|}{\textbf{L}} &
  \multicolumn{1}{l|}{\begin{tabular}[c]{@{}l@{}}Makes confusions with\\ the ``not'' operator.\end{tabular}} &
  \multicolumn{1}{l|}{\begin{tabular}[c]{@{}l@{}}Captures the meaning of\\ the ``not'' operator, but\\ needed help/made mistakes.\end{tabular}} &
  \begin{tabular}[c]{@{}l@{}}Applies all logical\\ connectives correctly\\ and without help.\end{tabular} \\ \hline
\multicolumn{1}{|l|}{\textbf{HLI}} &
  \multicolumn{1}{l|}{\begin{tabular}[c]{@{}l@{}} ~\\ Does not split the\\ proof in two parts.\\~ \end{tabular}} &
  \multicolumn{1}{l|}{\begin{tabular}[c]{@{}l@{}}Splits the proof in two parts,\\ maybe with a bit help. \\~ \end{tabular}} &
  \begin{tabular}[c]{@{}l@{}}Sees how to use\\ the ``and'' in the proof.\end{tabular} \\ \hline
\multicolumn{1}{|l|}{\textbf{S}} &
  \multicolumn{1}{l|}{\begin{tabular}[c]{@{}l@{}}Bad approach of proof,\\ little to no structure\\ visible.\end{tabular}} &
  \multicolumn{1}{l|}{\begin{tabular}[c]{@{}l@{}}Basic structure detected,\\ maybe some parts not\\ formulated well enough.\end{tabular}} &
  \begin{tabular}[c]{@{}l@{}}CLean and confident\\ structure.\end{tabular} \\ \hline
\multicolumn{1}{|l|}{\textbf{UE}} &
  \multicolumn{1}{l|}{\begin{tabular}[c]{@{}l@{}}Never used an\\ example.\end{tabular}} &
  \multicolumn{1}{l|}{Talks about examples.} &
  \begin{tabular}[c]{@{}l@{}}Uses examples\\ in a helpful way.\end{tabular} \\ \hline
\end{tabular}
}
\caption{Grading for Question 1.}
\label{tab:grading1}
}
\hfill
\parbox{.5\linewidth}{
\centering
\resizebox{.5\textwidth}{!}{
\begin{tabular}{|llll|}
\hline
\multicolumn{4}{|l|}{\cellcolor[HTML]{FFFF00}\textbf{Question 2}} \\ \hline
\multicolumn{1}{|l|}{} &
  \multicolumn{1}{l|}{0 Points} &
  \multicolumn{1}{l|}{1 Point} &
  2 Points \\ \hline
\multicolumn{1}{|l|}{\textbf{D}} &
  \multicolumn{1}{l|}{\begin{tabular}[c]{@{}l@{}}Needs help with the\\ definitions, even when\\ given.\end{tabular}} &
  \multicolumn{1}{l|}{\begin{tabular}[c]{@{}l@{}}Understands the\\ needed definitions, i.e.\\ some axioms, but maybe\\ does not know them all.\end{tabular}} &
  \begin{tabular}[c]{@{}l@{}}Uses them in a\\ confident way.\end{tabular} \\ \hline
\multicolumn{1}{|l|}{\textbf{MS}} &
  \multicolumn{1}{l|}{\begin{tabular}[c]{@{}l@{}}Struggles with writing\\ mathematics.\end{tabular}} &
  \multicolumn{1}{l|}{\begin{tabular}[c]{@{}l@{}}Has a cLean way of\\ writing mathematics.\end{tabular}} &
  \begin{tabular}[c]{@{}l@{}}No neglection\\ detected.\end{tabular} \\ \hline
\multicolumn{1}{|l|}{\textbf{L}} &
  \multicolumn{1}{l|}{\begin{tabular}[c]{@{}l@{}}Does not capture the\\ meaning of ``for all’.\end{tabular}} &
  \multicolumn{1}{l|}{\begin{tabular}[c]{@{}l@{}}Captures the meaning of\\ ``for all'', but\\ needed help/made mistakes.\end{tabular}} &
  \begin{tabular}[c]{@{}l@{}}Applies all logical\\ connectives correctly\\ and without help.\end{tabular} \\ \hline
\multicolumn{1}{|l|}{\textbf{HLI}} &
  \multicolumn{1}{l|}{\begin{tabular}[c]{@{}l@{}}Needs to be told to\\ use induction and\\ did not manage to\\ apply mul\_succ.\end{tabular}} &
  \multicolumn{1}{l|}{\begin{tabular}[c]{@{}l@{}}Understands to use\\ induction, maybe\\ with a bit help.\end{tabular}} &
  \begin{tabular}[c]{@{}l@{}}Sees how to use\\ mul\_succ.\end{tabular} \\ \hline
\multicolumn{1}{|l|}{\textbf{S}} &
  \multicolumn{1}{l|}{\begin{tabular}[c]{@{}l@{}}Bad approach of proof,\\ little to no structure\\ visible.\end{tabular}} &
  \multicolumn{1}{l|}{\begin{tabular}[c]{@{}l@{}}Basic structure detected,\\ maybe some parts not\\ formulated well enough.\end{tabular}} &
  \begin{tabular}[c]{@{}l@{}}CLean and confident\\ structure.\end{tabular} \\ \hline
\multicolumn{1}{|l|}{\textbf{UE}} &
  \multicolumn{1}{l|}{\begin{tabular}[c]{@{}l@{}}Never used an\\ example.\end{tabular}} &
  \multicolumn{1}{l|}{Talks about examples.} &
  \begin{tabular}[c]{@{}l@{}}Uses examples\\ in a helpful way.\end{tabular} \\ \hline
\end{tabular}
}
\caption{Grading for Question 2.}
\label{tab:grading2}
}
\end{table}

\begin{table}[H]
\parbox{.5\linewidth}{
\centering
\resizebox{.5\textwidth}{!}{
\begin{tabular}{|llll|}
\hline
\multicolumn{4}{|l|}{\cellcolor[HTML]{B3E5A1}\textbf{Question 3}}                                \\ \hline
\multicolumn{1}{|l|}{} & \multicolumn{1}{l|}{0 Points} & \multicolumn{1}{l|}{1 Point} & 2 Points \\ \hline
\multicolumn{1}{|l|}{\textbf{D}} &
  \multicolumn{1}{l|}{\begin{tabular}[c]{@{}l@{}}Needs help with the\\ definitions, even when\\ given.\end{tabular}} &
  \multicolumn{1}{l|}{\begin{tabular}[c]{@{}l@{}}Understands the\\ needed definitions, i.e.\\ refl, symm, trans, but maybe\\ does not know them all.\end{tabular}} &
  \begin{tabular}[c]{@{}l@{}}Uses them in a\\ confident way.\end{tabular} \\ \hline
\multicolumn{1}{|l|}{\textbf{MS}} &
  \multicolumn{1}{l|}{\begin{tabular}[c]{@{}l@{}}Struggles with writing\\ mathematics.\end{tabular}} &
  \multicolumn{1}{l|}{\begin{tabular}[c]{@{}l@{}}Has a cLean way of\\ writing mathematics.\end{tabular}} &
  \begin{tabular}[c]{@{}l@{}}No neglection\\ detected.\end{tabular} \\ \hline
\multicolumn{1}{|l|}{\textbf{L}} &
  \multicolumn{1}{l|}{\begin{tabular}[c]{@{}l@{}}Struggles with the\\ ``if-and-only-ifs’. \\~ \end{tabular}} &
  \multicolumn{1}{l|}{\begin{tabular}[c]{@{}l@{}}Captures the meaning\\ of ``if-and-only-ifs’, but\\ needed help/made mistakes.\end{tabular}} &
  \begin{tabular}[c]{@{}l@{}}Applies all logical\\ connectives correctly\\ and without help.\end{tabular} \\ \hline
\multicolumn{1}{|l|}{\textbf{HLI}} &
  \multicolumn{1}{l|}{\begin{tabular}[c]{@{}l@{}}Does not know how to\\ prove symm and trans.\end{tabular}} &
  \multicolumn{1}{l|}{\begin{tabular}[c]{@{}l@{}}Understands how to prove\\ symm and trans, maybe\\ with a bit help.\end{tabular}} &
  \begin{tabular}[c]{@{}l@{}}Sees how to prove\\ symm and trans. \\~ \end{tabular} \\ \hline
\multicolumn{1}{|l|}{\textbf{S}} &
  \multicolumn{1}{l|}{\begin{tabular}[c]{@{}l@{}}Bad approach of proof,\\ little to no structure\\ visible.\end{tabular}} &
  \multicolumn{1}{l|}{\begin{tabular}[c]{@{}l@{}}Basic structure detected,\\ maybe some parts not\\ formulated well enough.\end{tabular}} &
  \begin{tabular}[c]{@{}l@{}}CLean and confident\\ structure.\end{tabular} \\ \hline
\multicolumn{1}{|l|}{\textbf{UE}} &
  \multicolumn{1}{l|}{\begin{tabular}[c]{@{}l@{}}Never used an\\ example.\end{tabular}} &
  \multicolumn{1}{l|}{Talks about examples.} &
  \begin{tabular}[c]{@{}l@{}}Uses examples\\ in a helpful way.\end{tabular} \\ \hline
\end{tabular}
}
\caption{Grading for Question 3.}
\label{tab:grading3}
}
\hfill
\parbox{.5\linewidth}{
\centering
\resizebox{.5\textwidth}{!}{
\begin{tabular}{|llll|}
\hline
\multicolumn{4}{|l|}{\cellcolor[HTML]{95DCF7}\textbf{Question 4}}                                \\ \hline
\multicolumn{1}{|l|}{} & \multicolumn{1}{l|}{0 Points} & \multicolumn{1}{l|}{1 Point} & 2 Points \\ \hline
\multicolumn{1}{|l|}{\textbf{D}} &
  \multicolumn{1}{l|}{\begin{tabular}[c]{@{}l@{}}Needs help with the\\ definitions, even when\\ given.\end{tabular}} &
  \multicolumn{1}{l|}{\begin{tabular}[c]{@{}l@{}}Understands the needed\\ definitions, i.e. union and\\ intersection, but maybe\\ does not know them all.\end{tabular}} &
  \begin{tabular}[c]{@{}l@{}}Uses them in a\\ confident way.\end{tabular} \\ \hline
\multicolumn{1}{|l|}{\textbf{MS}} &
  \multicolumn{1}{l|}{\begin{tabular}[c]{@{}l@{}}Struggles with writing\\ mathematics.\end{tabular}} &
  \multicolumn{1}{l|}{\begin{tabular}[c]{@{}l@{}}Has a cLean way of\\ writing mathematics.\end{tabular}} &
  \begin{tabular}[c]{@{}l@{}}No neglection\\ detected.\end{tabular} \\ \hline
\multicolumn{1}{|l|}{\textbf{L}} &
  \multicolumn{1}{l|}{\begin{tabular}[c]{@{}l@{}}Does not capture the\\ ``does-it-follow’\\ correctly.\end{tabular}} &
  \multicolumn{1}{l|}{\begin{tabular}[c]{@{}l@{}}Captures the meaning\\ of ``how-it-follows’, but\\ needed help/made mistakes.\end{tabular}} &
  \begin{tabular}[c]{@{}l@{}}Follows through\\ correctly and\\ without help.\end{tabular} \\ \hline
\multicolumn{1}{|l|}{\textbf{HLI}} &
  \multicolumn{1}{l|}{\begin{tabular}[c]{@{}l@{}}Does not set $x$ in\\ union over $J$.\end{tabular}} &
  \multicolumn{1}{l|}{\begin{tabular}[c]{@{}l@{}}Understands how to\\ prove the subset, maybe\\ with a bit help.\end{tabular}} &
  \begin{tabular}[c]{@{}l@{}}Sees how to use the\\ definition of union\\ and intersection.\end{tabular} \\ \hline
\multicolumn{1}{|l|}{\textbf{S}} &
  \multicolumn{1}{l|}{\begin{tabular}[c]{@{}l@{}}Bad approach of proof,\\ little to no structure\\ visible.\end{tabular}} &
  \multicolumn{1}{l|}{\begin{tabular}[c]{@{}l@{}}Basic structure detected,\\ maybe some parts not\\ formulated well enough.\end{tabular}} &
  \begin{tabular}[c]{@{}l@{}}CLean and confident\\ structure.\end{tabular} \\ \hline
\multicolumn{1}{|l|}{\textbf{UE}} &
  \multicolumn{1}{l|}{\begin{tabular}[c]{@{}l@{}}Never used an\\ example.\end{tabular}} &
  \multicolumn{1}{l|}{Talks about examples.} &
  \begin{tabular}[c]{@{}l@{}}Uses examples\\ in a helpful way.\end{tabular} \\ \hline
\end{tabular}
}
\caption{Grading for Question 4.}
\label{tab:grading4}
}
\end{table}

   These grading tables are used as support for marking the interviews. Note that after individual grading, we had intense discussions about the scores given. During these discussions, some criteria in the grading tables may have been stretched slightly. For example, this might occur if a student did not write down a definition correctly, but it was clear during the interview that they understood the definition well enough.

   The questionnaire contains questions that provide insight into how Lean impacted the students who worked with us during the semester. The answers were collected using a Google form. Students were asked to respond to the following questions.

\begin{itemize}
    \item Did Lean motivate you to spend more time on a proof than usual? That is, would you spend more time trying to solve an exercise in Lean or in the Natural Number Game than on paper?
    \item Did Lean improve your proving skills in general?
    \item Did Lean influence how well you felt prepared for the final exam?
    \item How often did you use Lean besides the meetings?
    \item Should students be taught mathematics with Lean in the future?
    \item Open question: Do you plan to continue working with Lean during your studies? Why or why not?
    \item Open question: Would you have done something differently in the way we held our meetings with you?
\end{itemize}

  The first few questions could be answered with `No/Never', `Rarely/Hardly', `Sometimes/A bit', or `Yes/Often'. Question 5 offers the options `No', `No opinion', `That depends', or `Yes'. All non-open questions include an optional comment section. This questionnaire helps us to better analyze the students' progress. If a Lean student does not perform well in the exam and/or the interview, we check how these students answered the questionnaire. If they answered Question 3 with `Never' or `Rarely,' we know that the issue is not the meetings themselves, but a lack of motivation to work with Lean at home. The results of this questionnaire can be found in Section \ref{subsec:questionnaire}.

\section{Results}\label{sec:results}

Based on \cite{thoma2021learning}, we first compare the performance of the students participating in the interviews during the first exercise sheet and the exam to ensure that Lean students were not already better students before the sessions began. Then, the results from the interviews are presented and discussed. Finally, we examine some individual students or compare Lean and Non-Lean students with similar performance, and present the significance of the results using an independent $t$-test and a Mann-Whitney $U$-test.

\subsection{Exercise Sheets and the Exam} \label{subsec:exerciseandexamresults}

In this section, we present the mean and median scores from the exercise sheets and the final exam of the course ``Foundations of Mathematics''. We do this to compare the performance between Lean and Non-Lean students. Unlike the results from the interviews, this data includes the scores from all Non-Lean students in the course.

\begin{table}[H]
\centering
\resizebox{\textwidth}{!}{
\begin{tabular}{|l|lllllll|}
\hline
                        & \multicolumn{1}{l|}{\textbf{Sheet 1}} & \multicolumn{1}{l|}{\textbf{Sheet 2}} & \multicolumn{1}{l|}{\textbf{Sheet 3}} & \multicolumn{1}{l|}{\textbf{Sheet 4}} & \multicolumn{1}{l|}{\textbf{Sheet 5}} & \multicolumn{1}{l|}{\textbf{Sheet 6}} & \textbf{Exam} \\ \hline
\textbf{Lean}           & \multicolumn{7}{l|}{}                                                                                                                                                                                                                                         \\ \hline
Median                  & \multicolumn{1}{l|}{19}               & \multicolumn{1}{l|}{16.5}             & \multicolumn{1}{l|}{20}               & \multicolumn{1}{l|}{20}               & \multicolumn{1}{l|}{18.5}             & \multicolumn{1}{l|}{17.3}             & 47            \\ \hline
Mean                 & \multicolumn{1}{l|}{18.7}             & \multicolumn{1}{l|}{17.3}             & \multicolumn{1}{l|}{19.8}             & \multicolumn{1}{l|}{19.9}             & \multicolumn{1}{l|}{19.1}             & \multicolumn{1}{l|}{17.6}             & 46            \\ \hline
\textbf{Non-Lean}       & \multicolumn{7}{l|}{}                                                                                                                                                                                                                                         \\ \hline
Median                  & \multicolumn{1}{l|}{19}               & \multicolumn{1}{l|}{18.5}             & \multicolumn{1}{l|}{19}               & \multicolumn{1}{l|}{18.5}             & \multicolumn{1}{l|}{17.5}             & \multicolumn{1}{l|}{18}               & 35            \\ \hline
Mean                 & \multicolumn{1}{l|}{18.2}             & \multicolumn{1}{l|}{17.4}             & \multicolumn{1}{l|}{18.7}             & \multicolumn{1}{l|}{17.7}             & \multicolumn{1}{l|}{16.2}             & \multicolumn{1}{l|}{17.2}             & 34.9          \\ \hline
\textbf{Total possible} & \multicolumn{1}{l|}{20}               & \multicolumn{1}{l|}{20}               & \multicolumn{1}{l|}{20}               & \multicolumn{1}{l|}{20}               & \multicolumn{1}{l|}{20}               & \multicolumn{1}{l|}{20}             & 60            \\ \hline
\end{tabular}}
\caption{Lean and Non-Lean students exercise sheets and exam performance.}
\label{tab:resultsHao}
\end{table}

We have access to the anonymized scores of all the students, with Lean students labeled as \textbf{L}. Since not all students submitted every exercise sheet, we exclude incomplete data and compute the mean and median for only those students who submitted the exercise sheets. This way, we consider around thirty-five Non-Lean and five Lean students. For the exam, we use the final score, not the exam grade, so we do not need to account for the grading key used. It is evident that while Lean and Non-Lean students performed fairly equally on the first exercise sheet, Lean students scored over 10 points higher on average than Non-Lean students in the exam.

\subsection{Interviews} \label{subsec:interviewresults}
We begin by stating the mean, lowest, and highest scores (points for all six criteria summed) for each question, without comparing Lean and Non-Lean students.

\begin{table}[H]
\centering
\begin{tabular}{|l|l|l|l|}
\hline
\textbf{Question} & \textbf{Mean Score} & \textbf{Lowest Score} & \textbf{Highest Score} \\ \hline
\textbf{1}        & 4.4                    & 0                     & 9                     \\ \hline
\textbf{2}        & 6.4                    & 2                     & 10                     \\ \hline
\textbf{3}        & 8.1                    & 4                     & 10                     \\ \hline
\textbf{4}        & 7.7                    & 3                     & 12                     \\ \hline
\end{tabular}
\caption{Mean, lowest and highest score for each exercise of the interview.}
\label{tab:scoresperexercise}
\end{table}

Keep in mind that, with six criteria, the highest possible score is 12 points per exercise. We will discuss the results for each question in more detail later. Here are the total scores from the interviews over all four questions:

\begin{table}[H]
\centering
\resizebox{\textwidth}{!}{
\begin{tabular}{|l|l|l|l|l|l|l|l|l|l|l|l|l|}
\hline
\textbf{}                                                                  & \textbf{Laurin} & \textbf{Nevia} & \textbf{Ladina} & \textbf{Niculin} & \textbf{Lavinia} & \textbf{Nuot} & \textbf{Linard} & \textbf{Nicola} & \textbf{Liun} & \textbf{\begin{tabular}[c]{@{}l@{}}Lean \\ Mean\end{tabular}} & \textbf{\begin{tabular}[c]{@{}l@{}}Non-Lean \\ Mean\end{tabular}} & \textbf{Mean} \\ \hline
\textbf{Definitions}                                                       & 4               & 3              & 8               & 4             & 6               & 2               & 7                 & 8               & 4              & 5.8                                                             & 4.25                                                                & 5.1      \\ \hline
\textbf{\begin{tabular}[c]{@{}l@{}}Mathematical\\ Symbols\end{tabular}}    & 3               & 6              & 7               & 7             & 8               & 4               & 7                 & 6               & 5              & 6                                                             & 5.75                                                                & 5.9      \\ \hline
\textbf{Logic}                                                             & 4               & 6              & 8               & 5             & 7               & 2               & 6                 & 7               & 4              & 5.8                                                             & 5                                                                 & 5.4      \\ \hline
\textbf{\begin{tabular}[c]{@{}l@{}}High Level \\ Idea\end{tabular}}        & 4               & 5              & 8               & 4             & 7               & 1               & 4                 & 6               & 2              & 5                                                             & 4                                                                 & 4.6               \\ \hline
\textbf{\begin{tabular}[c]{@{}l@{}}Structure of \\ the Proof\end{tabular}} & 4               & 4              & 7               & 3             & 6               & 2               & 4                 & 7               & 4              & 5                                                             & 4                                                                & 4.6      \\ \hline
\textbf{Use of Examples}                                                   & 0               & 0              & 1               & 2             & 2               & 0               & 2                 & 2               & 1              & 1.2                                                             & 1                                                                   & 1.1      \\ \hline
\textbf{Total}                                                             & 19              & 24             & 39              & 25            & 36              & 11              & 30                & 36              & 20             & 28.8                                                            & 24                                                               & 26.7     \\ \hline
\end{tabular}}
\caption{Total scores in each criterion for all students.}
\label{tab:totalscores}
\end{table}

\begin{table}[H]
\centering
\resizebox{\textwidth}{!}{
\begin{tabular}{|l|l|l|l|l|l|l|l|l|l|l|l|l|}
\hline
\textbf{General Mean}                                                  & \textbf{Laurin} & \textbf{Nevia} & \textbf{Ladina} & \textbf{Niculin} & \textbf{Lavinia} & \textbf{Nuot} & \textbf{Linard} & \textbf{Nicola} & \textbf{Liun} & \textbf{\begin{tabular}[c]{@{}l@{}}Lean\\ Mean\end{tabular}} & \textbf{\begin{tabular}[c]{@{}l@{}}Non-Lean\\ Mean\end{tabular}} & \textbf{Mean} \\ \hline
\textbf{Definitions}                                                       & 1               & 0.75           & 2               & 1          & 1.5               & 0.5            & 2                 & 1.5            & 1            & 1.45                                                            & 1.06                                                             & 1.3      \\ \hline
\textbf{\begin{tabular}[c]{@{}l@{}}Mathematical \\ Symbols\end{tabular}}   & 0.75             & 1.5           & 1.75               & 1.75          & 2            & 1            & 1.75               & 1.5               & 1.25           & 1.5                                                            & 1.44                                                             & 1.5      \\ \hline
\textbf{Logic}                                                             & 1            & 1.5            & 2            & 1.25          & 1.75            & 0.5               & 1.5              & 1.75            & 1           & 1.45                                                           & 1.25                                                              & 1.4      \\ \hline
\textbf{\begin{tabular}[c]{@{}l@{}}High Level \\ Idea\end{tabular}}        & 1            & 1.25           & 2            & 1          & 1.75            & 0.25             & 1              & 1.5             & 0.5           & 1.25                                                           & 1                                                              & 1.1            \\ \hline
\textbf{\begin{tabular}[c]{@{}l@{}}Structure of \\ the Proof\end{tabular}} & 1            & 1              & 1.75            & 0.75          & 1.5               & 0.5             & 1              & 1.75               & 1              & 1.25                                                           & 1                                                             & 1.1     \\ \hline
\textbf{Use of Examples}                                                   & 0               & 0              & 0.25               & 0.5          & 0.5             & 0            & 0.5               & 0.5             & 0.25              & 0.3                                                            & 0.25                                                               & 0.3      \\ \hline
\textbf{Total Mean}                                                      & 4.75            & 6              & 9.75            & 6.25             & 9            & 2.75            & 7.5              & 9             & 5           & 7.2                                                           & 6                                                             & 6.7      \\ \hline
\end{tabular}}
\caption{Mean scores in each criterion for all the students.}
\label{tab:averagescores}
\end{table}

In Table \ref{tab:totalscores} we see the total scores for each criterion and the total score for the interview. We can see that the mean score for Lean students is higher in each criterion. The largest difference in averages is for the \textit{Definitions} criterion, followed by the \textit{Logic} criterion. This may be due to the cLean and precise structure that Lean encourages.

Table \ref{tab:averagescores} shows the mean scores for each criterion across all four exercises. For example, Nevia scored an average of 1.5 points in the logic criterion over all four exercises. In the last row, we provide the total average points for each criterion, which is equivalent to considering the mean score for each question in the interview. Therefore, Tables \ref{tab:totalscores} and \ref{tab:averagescores} do not present different results; they simply represent the same data in different formats.

From these tables, one can already see that Lean students performed slightly better than Non-Lean students. We will now look at each question in detail to see if this is the case for every question.

\subsubsection*{Question 1} \label{subsec:question1}

\textit{Prove the following deMorgan law without a truth table:}
    $$\neg(p \vee q) = (\neg p)\wedge (\neg q)$$ 
\begin{table}[H]
\centering
\resizebox{\textwidth}{!}{
\begin{tabular}{|l|l|l|l|l|l|l|l|l|l|l|l|l|}
\hline
\cellcolor[HTML]{E49EDD}Question   1 & Laurin & Nevia & Ladina & Niculin & Lavinia & Nuot & Linard & Nicola & Liun & \textbf{\begin{tabular}[c]{@{}l@{}}Averarge \\ Lean\end{tabular}} & \textbf{\begin{tabular}[c]{@{}l@{}}Mean \\ Non-Lean\end{tabular}} & \textbf{Mean} \\ \hline
D                                    & 0      & 0     & 2      & 0    & 0      & 0      & 2        & 2      & 1     & 1                                                                 & 0.5                                                                    & 0.8       \\ \hline
MS                                   & 0      & 2     & 2      & 2    & 2      & 0      & 1        & 1      & 1     & 1.2                                                               & 1.25                                                                 & 1.2       \\ \hline
L                                    & 0      & 1     & 2      & 1    & 1      & 0      & 1        & 2      & 1     & 1                                                                 & 1                                                                  & 1       \\ \hline
HLI                                  & 0      & 1     & 2      & 0    & 1      & 0      & 1        & 1      & 1     & 1                                                                 & 0.5                                                                  & 0.8       \\ \hline
S                                    & 0      & 0     & 1      & 0    & 0      & 0      & 1        & 2      & 0     & 0.4                                                               & 0.5                                                                  & 0.4       \\ \hline
UE                                   & 0      & 0     & 0      & 1    & 0      & 0      & 1        & 0      & 0     & 0.2                                                               & 0.25                                                                 & 0.2       \\ \hline
\textbf{Total}                       & 0      & 4     & 9      & 4    & 4      & 0      & 7        & 8      & 4     & 4.8                                                               & 4                                                                    & 4.4                \\ \hline
\end{tabular}}
\caption{Results from Question 1.}
\label{tab:resultsquestion1}
\end{table}

\begin{figure}[H]
\centering
    \subfigure[Nicola's solution.]{\includegraphics[width=0.45\textwidth]{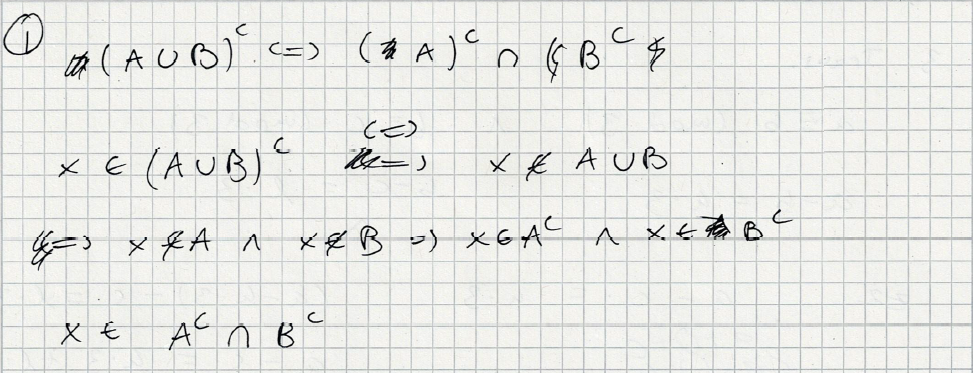}}
    \hspace{0.05\textwidth}
    \subfigure[Ladina's solution.]{\includegraphics[width=0.45\textwidth]{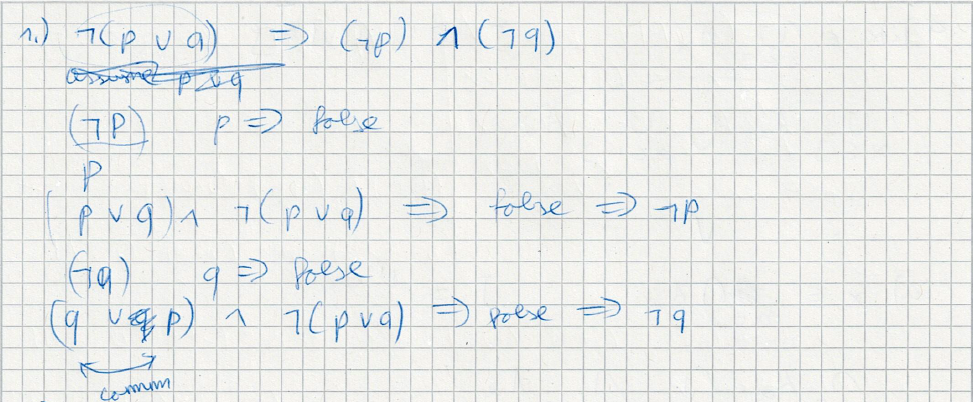}}
\caption{Two solutions for Question 1. Once with sets and once with propositional logic.}
\label{fig:question1selinaflurin}
\end{figure}

While Nicola manages to get eight points, he solves Question 1 using sets, Figure \ref{fig:question1selinaflurin}. He is not the only one to try that approach—Linard also uses sets to prove the statement. Though it is not the same setting, one could argue that, by the \textit{PAT interpretation}, an element of a set corresponds to a proof for a proposition, meaning De Morgan's law can be proven equivalently using set theory. Ladina is the only one able to solve this in the intended way. See Figure \ref{fig:question1selinaflurin} for Ladina's and Nicola's solutions. It turns out that solving Question 1 without a truth table is too difficult. Many make the mistake of trying to use the statement to prove itself, which occurs in a few of the interviews. Specifically, two Non-Lean students out of four and one Lean student out of five make this mistake. Interestingly, the Lean student is one of the lower performers, while the two Non-Lean students who make the mistake are average performers. Although Lean helps with understanding concepts like $\neg P \iff (P \Rightarrow \text{False})$, proving with first-order logic is an unfamiliar notation that takes time to master. It seems that more time with Lean is needed to develop proficiency in propositional logic proofs.

\subsubsection*{Question 2} \label{subsec:question2}

\textit{Using all axioms of addition and multiplication of natural numbers, prove:}
    $$\forall a, b, c \in \mathbb{N}, ~ a \cdot (b + c) = a \cdot b + a \cdot c $$

\begin{table}[H]
\centering
\resizebox{\textwidth}{!}{
\begin{tabular}{|l|l|l|l|l|l|l|l|l|l|l|l|l|}
\hline
\cellcolor[HTML]{FFFF66}Question   2 & Laurin & Nevia & Ladina & Niculin & Lavinia & Nuot & Linard & Nicola & Liun & \textbf{\begin{tabular}[c]{@{}l@{}}Average\\ Lean\end{tabular}} & \textbf{\begin{tabular}[c]{@{}l@{}}Mean\\ Non-Lean\end{tabular}} & \textbf{Mean} \\ \hline
D                                    & 1      & 0     & 2      & 0    & 2      & 0      & 1        & 2      & 1     & 1.4                                                              & 0.5                                                                 & 1                \\ \hline
MS                                   & 1      & 2     & 1      & 1    & 2      & 1      & 2        & 1      & 2     & 1.6                                                              & 1.25                                                                & 1.4       \\ \hline
L                                    & 1      & 1     & 2      & 1    & 2      & 1      & 2        & 2      & 2     & 1.8                                                              & 1.25                                                                & 1.6       \\ \hline
HLI                                  & 1      & 1     & 2      & 1    & 2      & 0      & 1        & 2      & 1     & 1.4                                                              & 1                                                                   & 1.2       \\ \hline
S                                    & 2      & 1     & 2      & 0    & 2      & 0      & 1        & 1      & 2     & 1.8                                                              & 0.5                                                                 & 1.2       \\ \hline
UE                                   & 0      & 0     & 0      & 0    & 0      & 0      & 0        & 0      & 0     & 0                                                                & 0                                                                   & 0                \\ \hline
\textbf{Total}                       & 6      & 5     & 9      & 3    & 10     & 2      & 7        & 8      & 8     & 8                                                                & 4.5                                                                 & 6.4       \\ \hline
\end{tabular}}
\caption{Results from Question 2.}
\label{tab:resultsquestion2}
\end{table}

While most students are able to figure out that a proof by induction is needed, some struggle with choosing the correct variable to perform induction on. For instance, in the following figure, Niculin attempts to do induction on all three variables.

\begin{figure}[H]
    \centering
    \includegraphics[width=0.6\textwidth]{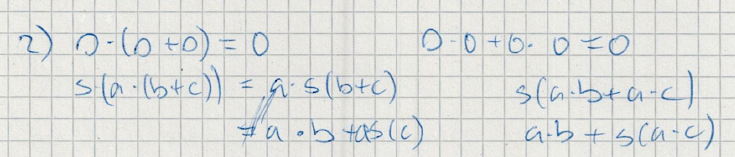}
    \caption{Niculin's solution for Question 2.}
    \label{fig:question2gian}
\end{figure}

  We check if they use the propositional equality correctly. As shown in Figures \ref{fig:question2laurinarmon} and \ref{fig:question2selinamadlaina} below, while both Lean and Non-Lean students both use the notion of propositional equality correctly, Nevia struggles more with the definition of successor multiplication than Lavinia does. Having played the Natural Number Game \cite{ref_NNG} may have given Lavinia an advantage.

\begin{figure}[H]
\centering
    \subfigure[Lavinia's solution.]{\includegraphics[width=0.45\textwidth]{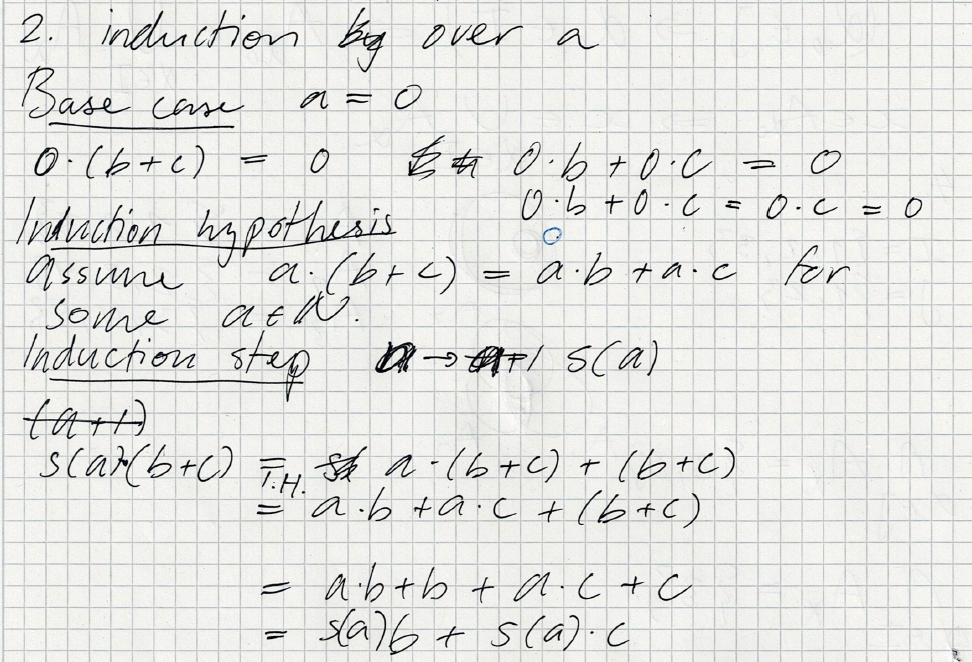}}
    \hspace{0.05\textwidth}
    \subfigure[Nevia's solution.]{\includegraphics[width=0.45\textwidth]{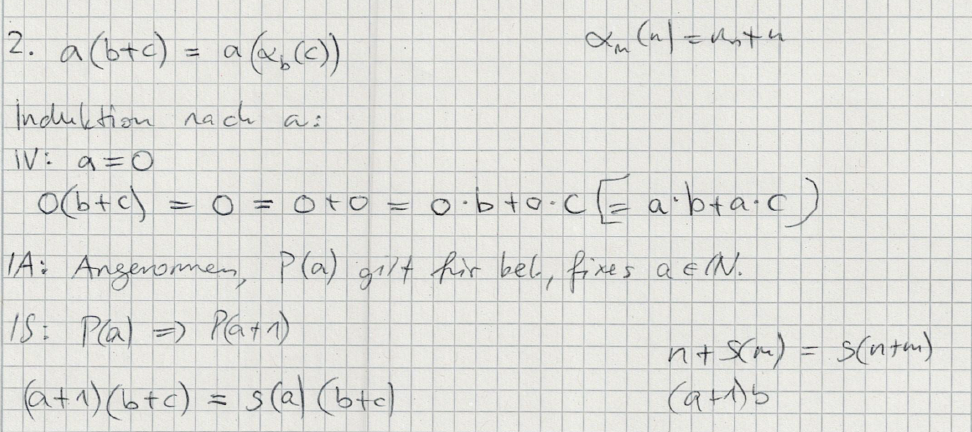}}
\caption{The propositional equality used correctly.}
\label{fig:question2laurinarmon}
\end{figure}

  In Figure \ref{fig:question2selinamadlaina} below, we see two students who incorrectly use propositional equality. It could be a simple notation mistake, but such errors can lead to difficulties in understanding the logic behind proving a statement.

\begin{figure}[H]
\centering
    \subfigure[Nicola's solution.]{\includegraphics[width=0.45\textwidth]{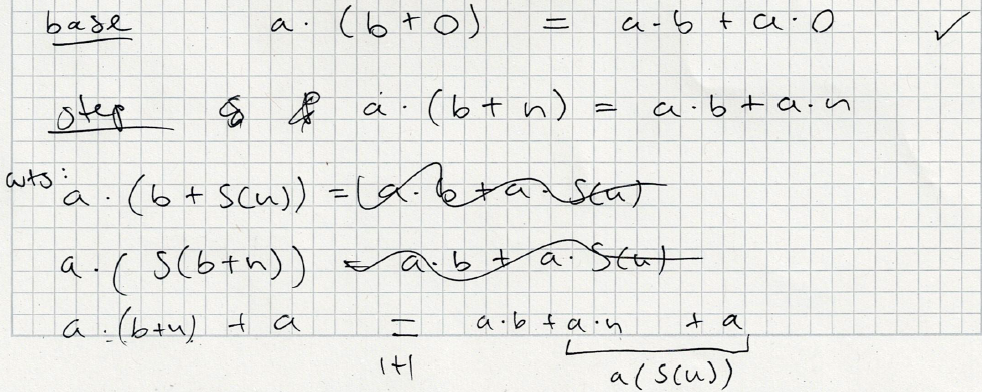}}
    \hspace{0.05\textwidth}
    \subfigure[Linard's solution.]{\includegraphics[width=0.45\textwidth]{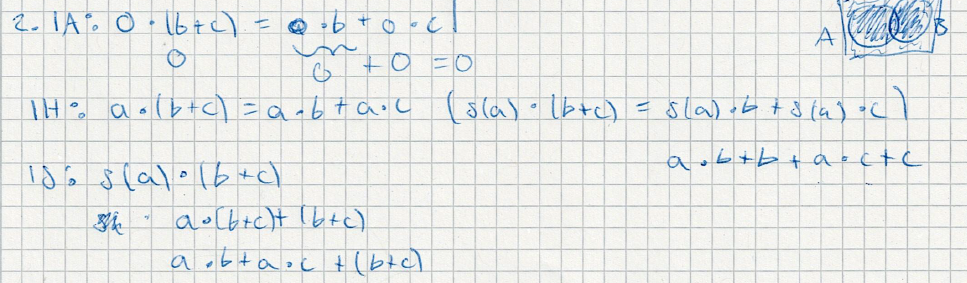}}
\caption{The propositional equality used incorrectly.}
\label{fig:question2selinamadlaina}
\end{figure}

\subsubsection*{Question 3} \label{subsec:question3}

\textit{Recall that: }
    $$a \equiv b \mod n \Leftrightarrow (\exists c, d, r \in \mathbb{N}, ~ a = c \cdot n + r  ~ \mbox{and} ~ b = d \cdot n + r) \Leftrightarrow (\exists k \in \mathbb{Z}, a - b = k \cdot n).$$
   \textit{ Show that} congruence mod $n$ \textit{with} $n = 3$ \textit{is an equivalence relation.}\\

\begin{table}[H]
\centering
\resizebox{\textwidth}{!}{
\begin{tabular}{|l|l|l|l|l|l|l|l|l|l|l|l|l|}
\hline
\cellcolor[HTML]{B5E6A2}Question   3 & Laurin & Nevia & Ladina & Niculin & Lavinia & Nuot & Linard & Nicola & Liun & \textbf{\begin{tabular}[c]{@{}l@{}}Averarge \\ Lean\end{tabular}} & \textbf{\begin{tabular}[c]{@{}l@{}}Mean \\ Non-Lean\end{tabular}} & \textbf{Mean} \\ \hline
D                                    & 2      & 2     & 2      & 2    & 2      & 1      & 2        & 2      & 2     & 2                                                                 & 1.75                                                                 & 1.9       \\ \hline
MS                                   & 2      & 1     & 2      & 2    & 2      & 2      & 2        & 2      & 1     & 1.8                                                               & 1.75                                                                 & 1.8       \\ \hline
L                                    & 2      & 2     & 2      & 2    & 2      & 0      & 1        & 2      & 0     & 1.4                                                               & 1.5                                                                  & 1.4       \\ \hline
HLI                                  & 2      & 1     & 2      & 1    & 2      & 0      & 1        & 2      & 0     & 1.4                                                               & 1                                                                    & 1.2       \\ \hline
S                                    & 2      & 2     & 2      & 2    & 2      & 1      & 1        & 2      & 2     & 1.8                                                               & 1.75                                                                 & 1.8       \\ \hline
UE                                   & 0      & 0     & 0      & 0    & 0      & 0      & 0        & 0      & 0     & 0                                                                 & 0                                                                    & 0                \\ \hline
\textbf{Total}                       & 10     & 8     & 10     & 9    & 10     & 4      & 7        & 10     & 5     & 8.4                                                               & 7.75                                                                 & 8.1       \\ \hline
\end{tabular}}
\caption{Results from Question 3.}
\label{tab:resultsquestion3}
\end{table}

\begin{figure}[H]
\centering
    \subfigure[Liun's solution.]{\includegraphics[width=0.45\textwidth]{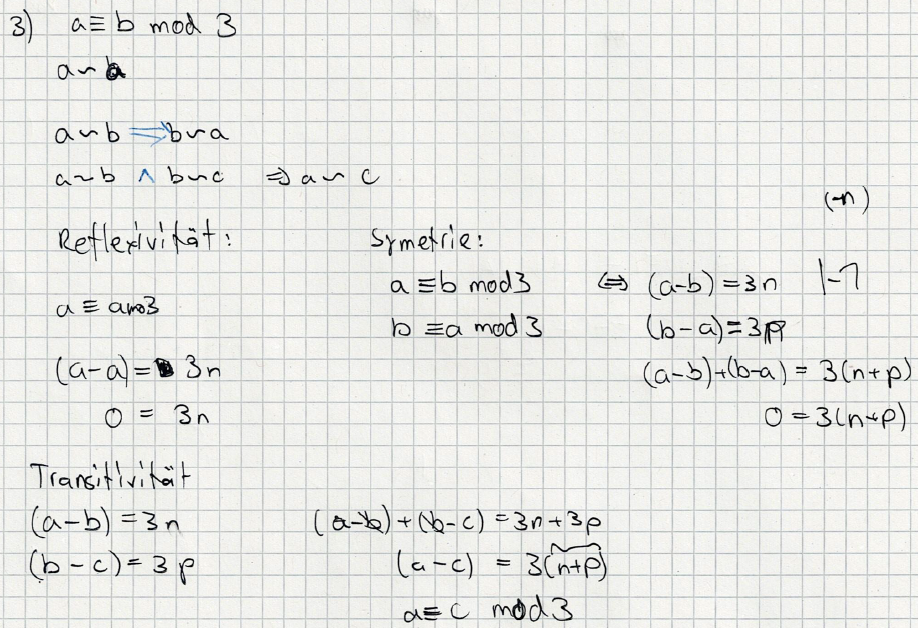}}
    \hspace{0.05\textwidth}
    \subfigure[Nuot's solution.]{\includegraphics[width=0.45\textwidth]{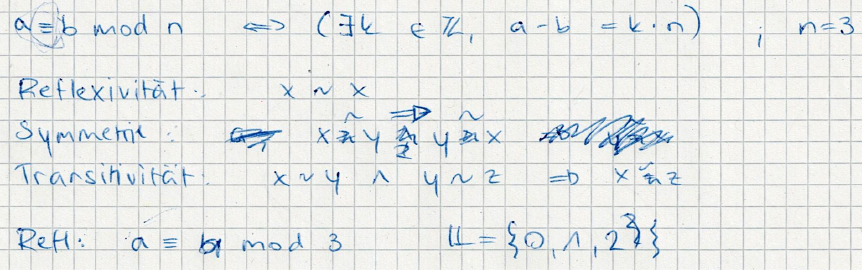}}
\caption{These two students struggle a bit with proving this statement.}
\label{fig:question3ursinlinard}
\end{figure}

Apart from a few students, everyone can correctly state the criteria needed to show that something is an equivalence relation. We compare two Lean and Non-Lean solutions for an adequate and a good proof of this statement.
In Figure \ref{fig:question3ursinlinard}, we seet that Nuot more or less knows the criteria for an equivalence proof, but he is unable to prove them. While Liun's proof is a bit more elaborate, he fails to prove symmetry correctly, as he tries to use the same method as for transitivity.

\begin{figure}[H]
\centering
    \subfigure[Niculin's solution.]{\includegraphics[width=0.45\textwidth]{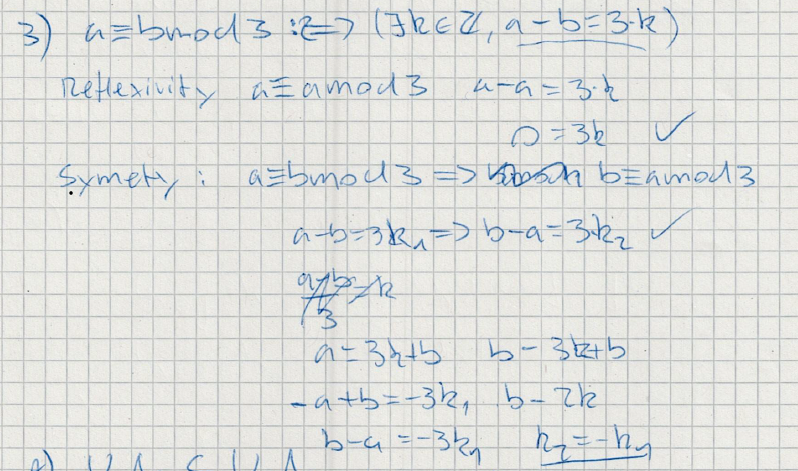}}
    \hspace{0.05\textwidth}
    \subfigure[Laurin's solution.]{\includegraphics[width=0.45\textwidth]{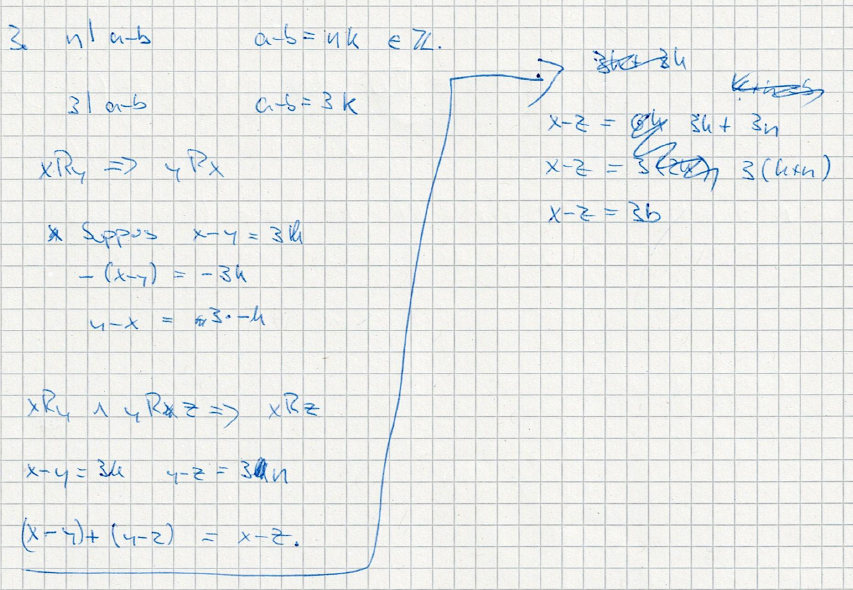}}
\caption{Two good examples of a proof for Question 3.}
\label{fig:question3gianandrin}
\end{figure}

  Niculin's and Laurin's solutions are nice examples of the proof in Figure \ref{fig:question3gianandrin}, but like other students too, they struggle to choose two different integers in the symmetry proof. They have it written in their proofs here, but they needed help to understand that they should use different integers in the proof.

\subsubsection*{Question 4} \label{subsec:question4}

\textit{If} $J \neq \emptyset$ \textit{and} $J \subseteq I$\textit{, does it follow that} $\bigcup_{\alpha \in J} A_{\alpha} \subseteq \bigcup_{\alpha \in I} A_{\alpha}$\textit{? What about} $\bigcap_{\alpha \in J} A_{\alpha} \subseteq \bigcap_{\alpha \in I} A_{\alpha}$\textit{?}\\

\begin{table}[H]
\centering
\resizebox{\textwidth}{!}{
\begin{tabular}{|l|l|l|l|l|l|l|l|l|l|l|l|l|}
\hline
\cellcolor[HTML]{94DCF8}Question   4 & Laurin & Nevia & Ladina & Niculin & Lavinia & Nuot & Linard & Nicola & Liun & \textbf{\begin{tabular}[c]{@{}l@{}}Averarge \\ Lean\end{tabular}} & \textbf{\begin{tabular}[c]{@{}l@{}}Mean \\ Non-Lean\end{tabular}} & \textbf{Mean} \\ \hline
D                                    & 1      & 1     & 2      & 2    & 2      & 1      & 2        & 2      & 0     & 1.4                                                               & 1.5                                                                  & 1.4       \\ \hline
MS                                   & 0      & 1     & 2      & 2    & 2      & 1      & 2        & 2      & 1     & 1.4                                                               & 1.5                                                                  & 1.4       \\ \hline
L                                    & 1      & 2     & 2      & 1    & 2      & 1      & 2        & 1      & 1     & 1.6                                                               & 1.25                                                                 & 1.4       \\ \hline
HLI                                  & 1      & 2     & 2      & 2    & 2      & 1      & 1        & 1      & 0     & 1.2                                                               & 1.5                                                                  & 1.3       \\ \hline
S                                    & 0      & 1     & 2      & 1    & 2      & 1      & 1        & 2      & 0     & 1                                                                 & 1.25                                                                 & 1.1       \\ \hline
UE                                   & 0      & 0     & 1      & 1    & 2      & 0      & 1        & 2      & 1     & 1                                                                 & 0.75                                                                 & 0.9       \\ \hline
\textbf{Total}                       & 3      & 7     & 11     & 9    & 12     & 5      & 9        & 10     & 3     & 7.6                                                               & 7.75                                                                 & 7.7       \\ \hline
\end{tabular}}
\caption{Results from Question 4.}
\label{tab:resultsquestion4}
\end{table}

  This is the only question in which Non-Lean students perform slightly better than Lean students. As we stated in Section \ref{subsec:interviews}, we expected that Non-Lean students have an advantage over Lean students in this question.

  Below we see Lavinia's approach to the proof. It is the highest-scoring question in all the interviews, as she uses a very nice example to disprove the statement about the intersection.

\begin{figure}[H]
\centering
    \subfigure[Proof for the union.]{\includegraphics[width=0.45\textwidth]{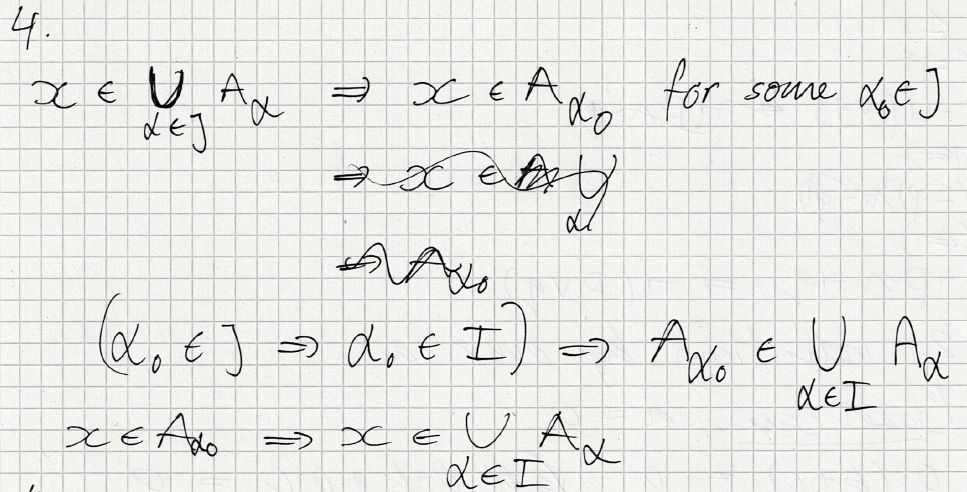}}
    \hspace{0.05\textwidth}
    \subfigure[Proof for the intersection.]{\includegraphics[width=0.45\textwidth]{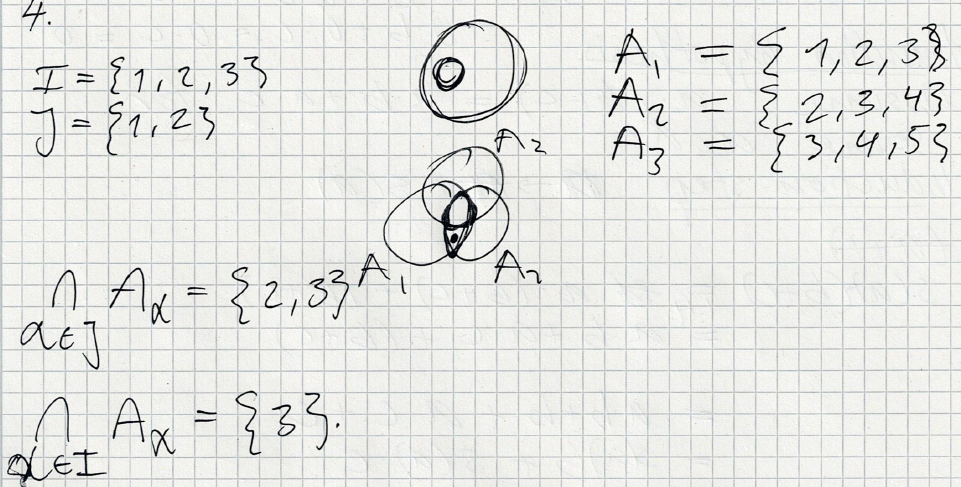}}
\caption{Lavinia's proof to Question 4.}
\label{fig:question4laurin}
\end{figure}

While in the figure above we see that Lavinia confidently proves the statement about the union and has no problem defining where the chosen $x$ should lie, for other students, both Lean and Non-Lean, this could lead to issues, as seen in the following figures.

\begin{figure}[H]
\centering
    \subfigure[Liun's solution.]{\includegraphics[width=0.45\textwidth]{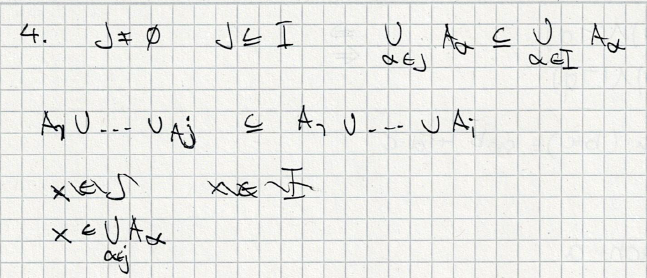}}
    \hspace{0.05\textwidth}
    \subfigure[Nuot's solution.]{\includegraphics[width=0.45\textwidth]{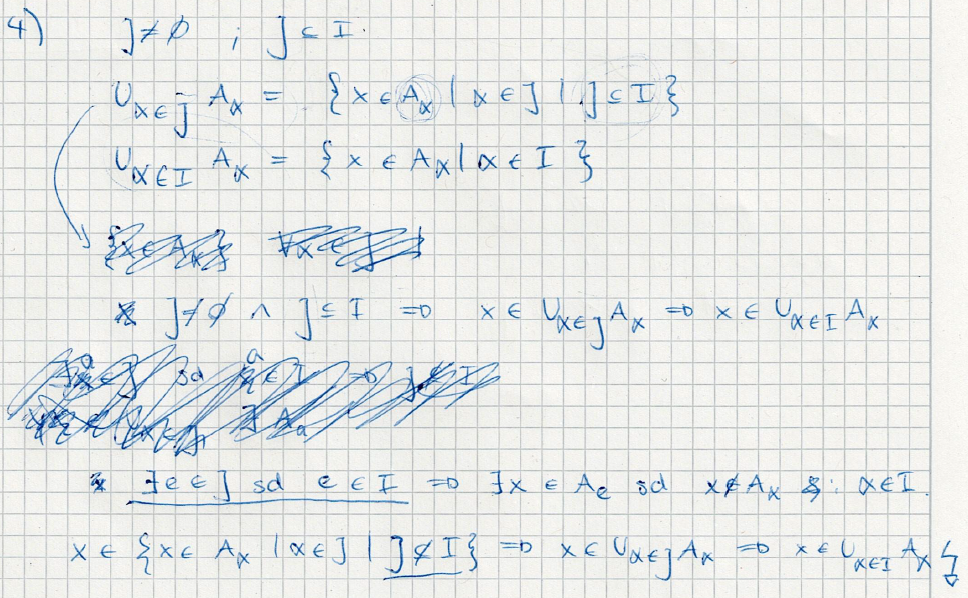}}
\caption{Both students struggle with finishing the proof.}
\label{fig:question4ursinlinard}
\end{figure}

The main reason Nuot is confused is that he thinks he remembers the first statement from the exercise sheets to be true, a common mistake that first-year students make—trying to remember the correctness of statements by heart. While he chooses a certain index $e$ to explain that $x$ would lie in some $A_\alpha$, he still struggles to define where $x$ lies exactly in the next line. Liun, on the other hand, cannot see the importance of choosing a specific label $\alpha_0$ to argue that $x$ would also lie in the union over the index set $I$. This is why he cannot finish the proof.

This problem of using the variable name given in the statement and not defining new labels is also seen in Question 3, where students get confused in the symmetry proof when they have to choose a different $k$. While we expect Lean students to avoid these kinds of mistakes, more time working with Lean is needed to fully eliminate them. As shown in the questionnaire, the more Lean students used Lean at home, the less they made mistakes like the ones mentioned here.

\subsection{Students' General Performance} \label{subsec:studentsgeneral}

   Having a small data set allows us to compare the students' performances more carefully. Sometimes we elaborate on individual students' performances, and sometimes we compare Lean and Non-Lean students with similar interview scores to highlight differences or similarities.

  We begin by discussing Ladina, as she is quite a special case. Ladina already has a lot of experience with Lean. She holds a Bachelor's in Informatics and has even implemented some formalizations of mathematics in Lean at some point. However, she started studying mathematics as a first-year student like the others. Coming from a computer science background, it is clear in various instances that she views mathematics differently than the other students and even us.

\begin{figure}[H]
    \centering
    \includegraphics[width=0.8\textwidth]{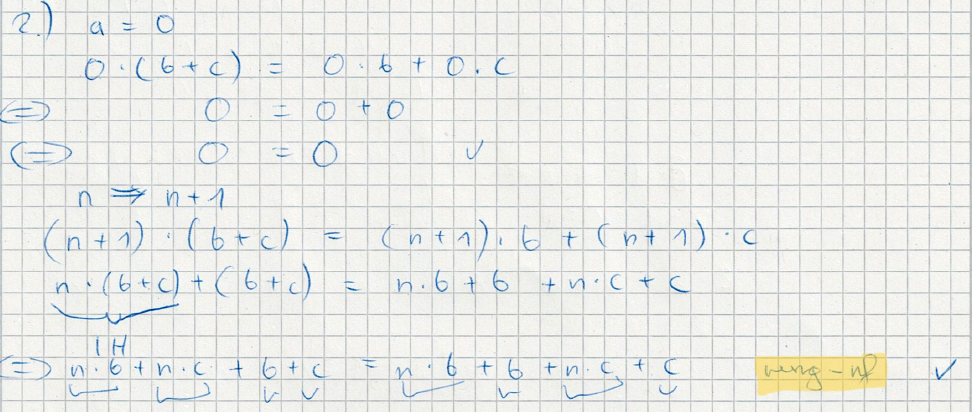}
    \caption{Question 2 from Ladina's interview.}
    \label{fig:question2flurin}
\end{figure}

   At the end of Question 2, Ladina mentions that she would simply use \textit{ring\_nf}, one of the Lean tactics we discussed earlier. This is not the only time she mentions a Lean tactic she would use to solve the next step of a proof.

 Although he is a Non-Lean student, Nicola scores better than most students in the interview. We do not know much about him, as he did not attend the meetings, but we learn at the beginning of the interview that he is a first-semester student. There are always students who are really good at mathematics, and it should be clear that the difference between Lean and Non-Lean students is not that all Lean students perform better than all Non-Lean students. Nicola has a clean and structured proving style without ever having used Lean. However, he solves the first exercise in the interview using sets, i.e., set theory, as seen in Figure \ref{fig:question1selinaflurin}. While he solves it correctly, it would have been much easier to do it using propositional logic, which he did not learn in the course ``Foundations of Mathematics''. In his example, we can see that no matter how good we are at something, learning something new is always beneficial.

Laurin and Nuot are the Lean and Non-Lean students, respectively, who score the lowest marks during the interviews. Based on feedback from Laurin, we understand that his score may be influenced by his phobia of oral examinations. We also see from the questionnaire that he is the one who answered `never' when asked if he used Lean outside of the meetings, which is why we hesitate to say that Lean did not improve his mathematical skills.

\begin{figure}[H]
\centering
    \subfigure[Laurin's solution.]{\includegraphics[width=0.35\textwidth]{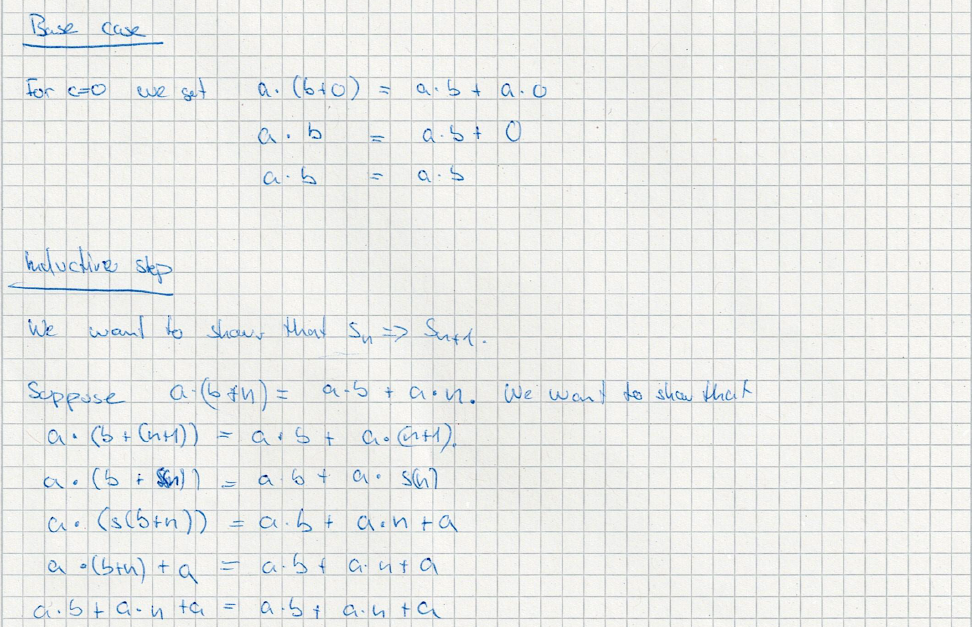}}
    \hspace{0.05\textwidth}
    \subfigure[Nuot's solution.]{\includegraphics[width=0.35\textwidth]{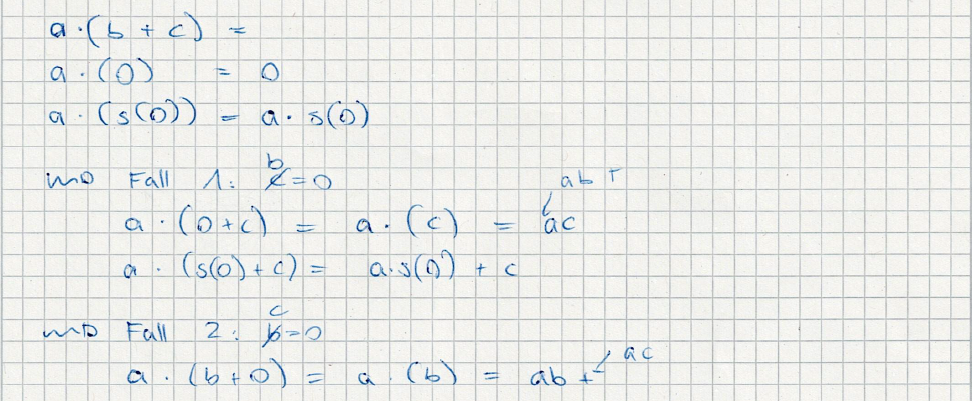}}
\caption{A Lean and Non-Lean solution to Question 2.}
\label{fig:question2andrinlindard}
\end{figure}

 Even though they both score lower in the interview than other students, one can see that Laurin writes his proof in Question 2 in a more structured way. This is possibly due to the Lean level covering natural induction.

Linard is one of the Lean students with a very high score and has the second-best score for \textit{Definitions} after Ladina. Interestingly, he also tries to solve Exercise 1 from the interview using set theory, as shown in Figure \ref{fig:question1madlaina} below. This suggests that first-order logic is not very prevalent among Lean students either. We would not attribute this to Lean being an incompetent theorem prover, but rather to the limited time spent learning about first-order logic in Lean.

\begin{figure}[H]
    \centering
    \includegraphics[width=0.9\textwidth]{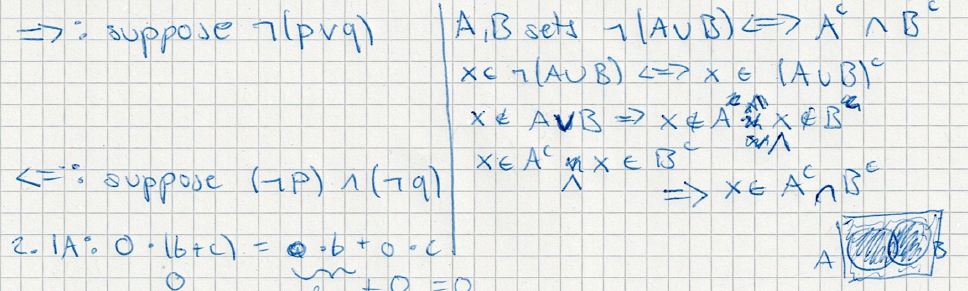}
    \caption{Question 1 from Linard's interview.}
    \label{fig:question1madlaina}
\end{figure}

Linard starts Question 1 with a propositional logic approach but later finishes it using sets. Note his cLean proof structure and the example he uses at the bottom right. It is one of the very few examples used across all interviews.

  Next, we would like to discuss Nevia, a Non-Lean student. Interestingly, her \textit{Logical thinking} is remarkable, but she struggles a bit with the \textit{Definitions} of mathematical objects, such as the union of indexed sets, Table \ref{tab:totalscores}. With some practice doing Lean exercises, she could certainly benefit and gain a better understanding of those definitions.

  Niculin is another Non-Lean student who scores almost the same mark as Nevia, but he seems much more confident. He spends more time on Exercises 1 and 3 than the other students because he gets a bit confused there. However, he is quite confident with definitions and proof structure throughout the entire interview.
  
Lavinia is a very enthusiastic Lean student. While she is not as experienced as Ladina, she did the most Lean practice at home compared to the other students. Her high score in the interview reflects her commitment. Keep in mind that she is the only student to score the full 12 points on a question. Even though she has never studied another course of study before, she has a strong interest in computers. One can clearly see some Lean-thinking during the interview, but unfortunately, she also struggles with Exercise 1. We would like to highlight how she solves Exercise 4. She uses an excellent example to disprove Subquestion 4.\ b), as shown in Figure \ref{fig:question4laurin}, and solves the exercise in a more elegant way than others. It is hard to say whether this is due to Lean or not. Besides Ladina, Lavinia is the most likely to continue using Lean regularly during her studies.

That leaves us with Liun, the last Lean student. He is a very committed student and learner, and he is the one who `praised' Lean the most in the questionnaire, mentioning that he often uses our Lean exercises to understand a proof he cannot figure out by himself. With a bit more experience, he could reach the level of Lavinia and Ladina. We are not sure whether Lean is more of a helpful tool for him or if his interest in proof assistants was also sparked by the sessions, but we hope he will continue using Lean either way.

If we compare the interviews overall, we see that Lean students perform slightly better, with outliers in both directions. Our sample is too small to draw any significant conclusions about the effect of teaching with Lean, but it is still encouraging to see that Lean students perform well in general. 

\subsection{Significance}
\label{sub:significance}

In this chapter, we present the tests we perform to check whether our results are significant. Our null hypothesis is that the mean performances of Lean and Non-Lean students are the same. We calculate the $p$-value to determine if we can reject the null hypothesis. We set the significance level at $\alpha = 0.05$.

  For the interview results and the final exam scores, we perform a $t$-test for two independent samples: Lean and Non-Lean students. We assume the variances to be heterogeneous, which affects the degrees of freedom.

   The $t$-value is computed as

$$t = \frac{\overline{x}_1-\overline{x}_2}{\sqrt{\frac{s_1^2}{n_1}+\frac{s_2^2}{n_2}}},$$~

   where $\overline{x}_k$ is the respective mean value for Lean or Non-Lean students, $s_k^2$ is the standard deviation and $n_k$ is the number of students per sample.

   The standard deviation is computed as follows,

$$s_k^2= \frac{1}{n_k-1}\sum_{i=1}^{n_k}(x_i-\overline{x}_k)^2,$$~ 

   where, $x_i$ stands for the performance of each individual student.

   To compute the $p$-value, we need the number of degrees of freedom. In our case, this number is computed as

$$df = \frac{(\frac{s_1^2}{n_1}+\frac{s_2^2}{n_2})^2}{\frac{1}{n_1-1}(\frac{s_1^2}{n_1})^2+\frac{1}{n_2-1}(\frac{s_2^2}{n_2})^2}.$$

   In Tables \ref{tab:ttestinterviews} and \ref{tab:ttestexams} below, we present the given data.

\begin{table}[H]
\centering
\begin{tabular}{|l|l|l|l|l|l|l|l|}
\hline
 \textbf{$\overline{x}_1$} & \textbf{$\overline{x}_2$} & \textbf{$n_1$} & \textbf{$n_2$} & \textbf{$s_1^2$} & \textbf{$s_2^2$} & \textbf{$t$} & \textbf{$df$} \\ \hline
28.8 & 24 & 5 & 4 & 82.7 & 104.67 & 0.73 & 6.15 \\ \hline
\end{tabular}
\caption{Quantities from the interviews to find $p$-value.}
\label{tab:ttestinterviews}
\end{table}

\begin{table}[H]
\centering
\begin{tabular}{|l|l|l|l|l|l|l|l|}
\hline
 $\overline{x}_1$ & \textbf{$\overline{x}_2$} & \textbf{$n_1$} & \textbf{$n_2$} & \textbf{$s_1^2$} & \textbf{$s_2^2$} & \textbf{$t$} & \textbf{$df$} \\ \hline
 46 & 34.92 & 5 & 52 & 32.5 & 181.17 & 3.51 & 9.23 \\ \hline
\end{tabular}
\caption{Quantities from the exams to find $p$-value.}
\label{tab:ttestexams}
\end{table}

   Using an online tool, we calculate the $p$-values. For this, we just have to type in the $t$-values and the degrees of freedom for both tests \cite{ref_tdistr}.\\

   The following $p$-values were found.

\begin{table}[H]
    \centering
    \begin{tabular}{|l|l|l|}
    \hline
    & \textbf{Interviews} & \textbf{Exams} \\ \hline
    \textbf{$p$-values} & 0.49 & 0.006 \\ \hline
    \end{tabular}
    \caption{$p$-values of interview and exam performances.}
    \label{tab:pvalues}
\end{table}

 We now observe that, for the interviews, while Lean students perform better on average than Non-Lean students, the difference in the means is not significant. For the exam, however, the $p$-value is smaller than the significance level. Therefore, we can reject the null hypothesis and state that Lean students performed significantly better on average in the exam than Non-Lean students.

  For the $t$-test, we assume the data to be normally distributed. To account for the possibility that the data may not be normally distributed, we conduct the nonparametric Mann-Whitney $U$-test for independent data \cite{ref_manntest}. This test involves several intermediate steps. These steps are explained next and are described for both the interview and exam results simultaneously.

 First, we rank all scores from $1$ to $n$, where $1$ corresponds to the lowest score and $n$ to the highest. If there are tied ranks, we take the mean of the ranks. For example, if the ranks $6$, $7$, and $8$ all have the same score, they are all assigned the rank $7$. We then compute the sum of the Lean and Non-Lean ranks, denoted $T_1$ and $T_2$, respectively, and calculate a value needed later, which we will call $tR$ (tied ranks).

$$tR = \sum_{i=1}^k \frac{t_i^3-t_i}{12}.$$

   Here, $k$ stands for the number of tied ranks and $t_i$ is the number of people sharing the same rank.

\begin{table}[H]
\centering
\begin{tabular}{|l|l|l|l|l|}
\hline
\textbf{$T_1$} & \textbf{$T_2$} & \textbf{$n_1$} & \textbf{$n_2$} & \textbf{$tR$} \\ \hline
27.5 & 17.5 & 5 & 4 & 0.5 \\ \hline
\end{tabular}
\caption{Rank sums and tied rank number from the interviews.}
\label{tab:ranksuminterview}
\end{table}

\begin{table}[H]
\centering
\begin{tabular}{|l|l|l|l|l|}
\hline
\textbf{$T_1$} & \textbf{$T_2$} & \textbf{$n_1$} & \textbf{$n_2$} & \textbf{$tR$} \\ \hline
217.5 & 1435.5 & 5 & 52 & 26.5 \\ \hline
\end{tabular}
\caption{Rank sums and tied rank number from the exams.}
\label{tab:ranksumexam}
\end{table}

   Next, we calculate the $U$-values for Lean and Non-Lean students and choose the smaller of the two as our $U$-value. We also compute $\mu_U$ and $\sigma_U$ as below

$$U_j = n_1 \cdot n_2 + \frac{n_j\cdot(n_j+1)}{2}-T_j,$$

$$\mu_U = \frac{n_1\cdot n_2}{2},$$

$$\sigma_U = \sqrt{\frac{n_1\cdot n_2}{n \cdot (n-1)}} \cdot \sqrt{\frac{n^3-n}{12}-tR}.$$~

 In these formulas, $j$ can take the values $1$ and $2$ for Lean and Non-Lean students, respectively, and $n$ is the sum of both sample sizes.

   The $z$-values in Tables \ref{tab:zdistrinterview} and \ref{tab:zdistrexam} are calculated using $z = \frac{U-\mu_U}{\sigma_U}$. Once we have that, we use an online calculator to find the $p$-value \cite{ref_zdistr}.

\begin{table}[H]
\centering
\begin{tabular}{|l|l|l|l|}
\hline
\textbf{$U$} & \textbf{$\mu_U$} & \textbf{$\sigma_U$} & \textbf{$z$}\\ \hline
7.5 & 10 & 4.07 & -0.61\\ \hline
\end{tabular}
\caption{Quantities from the interviews to find $p$-value.}
\label{tab:zdistrinterview}
\end{table}

\begin{table}[H]
\centering
\begin{tabular}{|l|l|l|l|}
\hline
\textbf{$U$} & \textbf{$\mu_U$} & \textbf{$\sigma_U$} & \textbf{$z$} \\ \hline
57.5 & 130 & 35.42 & -2.05 \\ \hline
\end{tabular}
\caption{Quantities from the exams to find $p$-value.}
\label{tab:zdistrexam}
\end{table}

 The Mann-Whitney $U$-test returns the same conclusions as the $t$-test. The better performance of Lean students during the interviews is not significant ($p > \alpha$), but for the exams, our results are significant ($p < \alpha$).

\begin{table}[H]
\centering
\begin{tabular}{|l|l|l|}
\hline
\textbf{} & \textbf{Interviews} & \textbf{Exams} \\ \hline
\textbf{$p$-values} & 0.54 & 0.04 \\ \hline
\end{tabular}
\caption{$p$-values of interview and exam performances.}
\label{tab:pvaluefromzdistr}
\end{table}

\subsection{Results from the Questionnaire} \label{subsec:questionnaire}

Finally, we present the results of the questionnaire from Chapter \ref{subsec:interviews}.

\subsubsection*{Did Lean motivate you to spend more time on a proof than usual? i.e.\ would you spend more time trying to solve an exercise in Lean (or the NNG) than on paper?}

\begin{figure}[H]
    \centering
    \includegraphics[width = 0.7\textwidth]{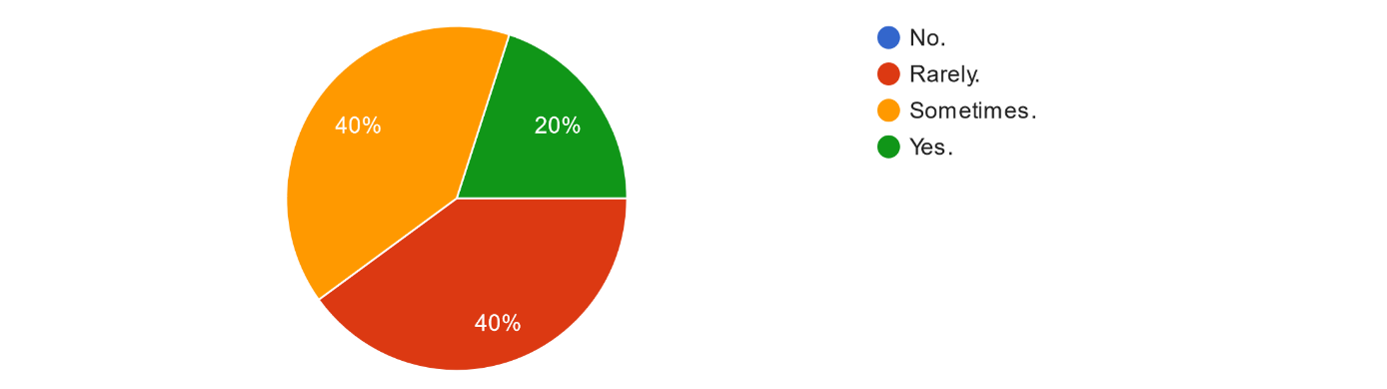}
    \caption{Question 1.}
    \label{fig:questionnaire1}
\end{figure}

   \textbf{Comments:} \begin{itemize}
    \item Lean kind of forces your hand in this case. However, occasionally it also happens by choice, as the exploration of the structure of the proof is more in-depth with Lean. This is not always the case, as some proofs pretty much can be auto-generated and it also depends on my interest in the proof.
    \item Especially the NNG.
\end{itemize}

\subsubsection*{Did Lean improve your proving skills in general?}

\begin{figure}[H]
    \centering
    \includegraphics[width = 0.7\textwidth]{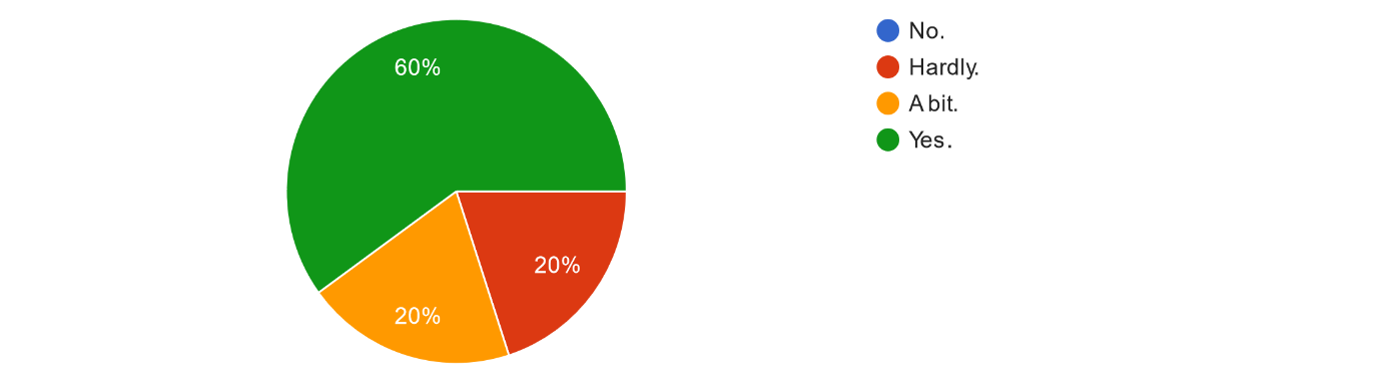}
    \caption{Question 2.}
    \label{fig:questionnaire2}
\end{figure}

   \textbf{Comments:} \begin{itemize}
    \item Absolutely, in typical proof patterns (inj/surj proofs, set equality, induction, contradiction, etc.) I've learned to start out in a very mechanical way, almost like mathematical muscle memory. I attribute this to Lean.
    \item Helped with structure of proofs.
\end{itemize}~

\subsubsection*{Did Lean have an influence on how well you felt prepared for the final exam?}

\begin{figure}[H]
    \centering
    \includegraphics[width = 0.8\textwidth]{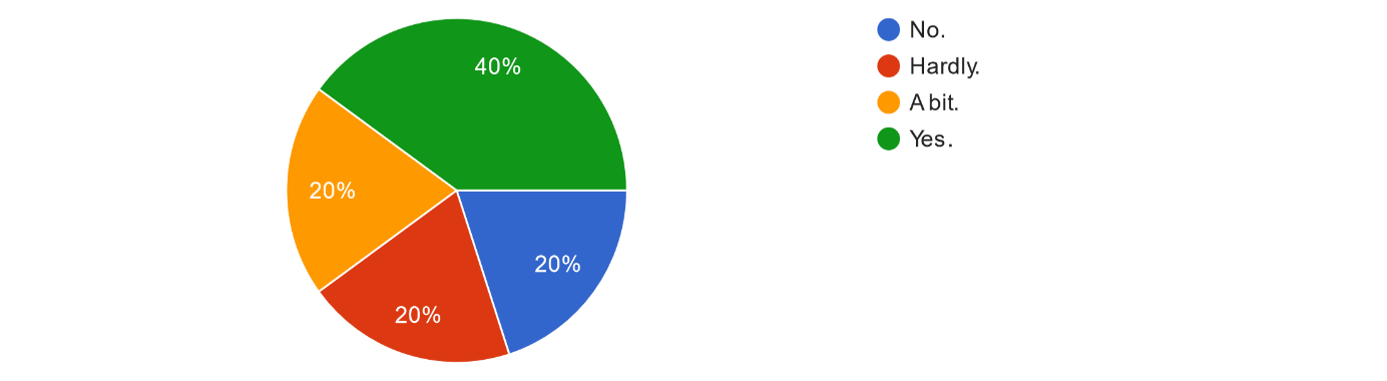}
    \caption{Question 3.}
    \label{fig:questionnaire3}
\end{figure}

   \textbf{Comments:} \begin{itemize}
    \item I assumed that I would be able to do the aforementioned mechanical proofs quite a bit faster than if I had not looked at Lean. However, in general, this was not the case and it showed in the exam. My time in the exam was mostly taken up with annoying algebraic manipulations and manipulations over summation symbols (stuff I use simp for). I blame this on the exam though, as the hard part should not be those aspects.
\end{itemize}

\subsubsection*{How often did you use Lean besides the meetings?}

\begin{figure}[H]
    \centering
    \includegraphics[width = 0.8\textwidth]{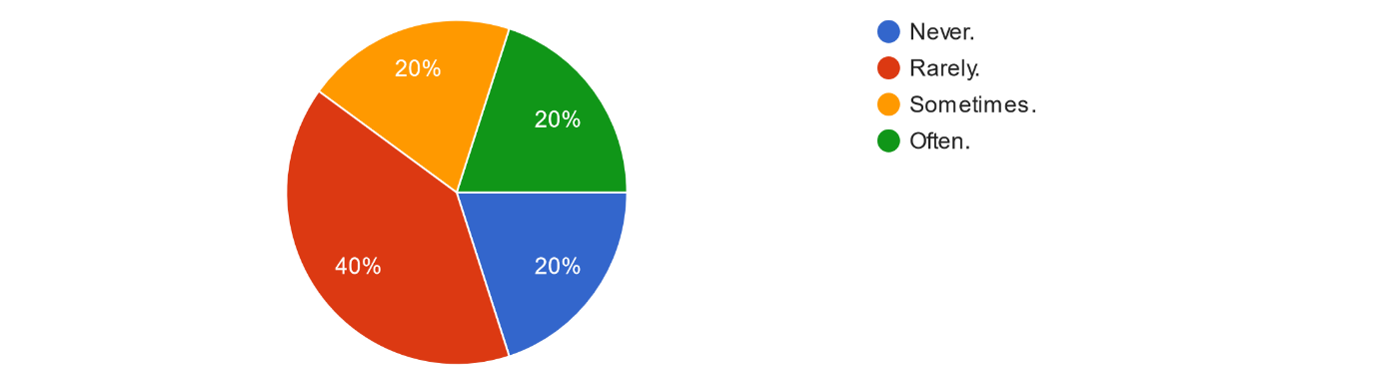}
    \caption{Question 4.}
    \label{fig:questionnaire4}
\end{figure}

   \textbf{Comments:} \begin{itemize}
    \item I try to formalize various theorems and am actively doing formal verification of code. Sometimes, when a theorem from one of my modules is not clear, I browse mathlib for it to see all the parameters it takes and the output, so I get a very structured overview of what the theorem exactly does (and I get to see applications of it in mathlib).
    \item Whenever I was stuck with an exercise I tried to compare with the code given by you.
\end{itemize}

\subsubsection*{Should students be taught mathematics with Lean in the future?}

\begin{figure}[H]
    \centering
    \includegraphics[width = 0.8\textwidth]{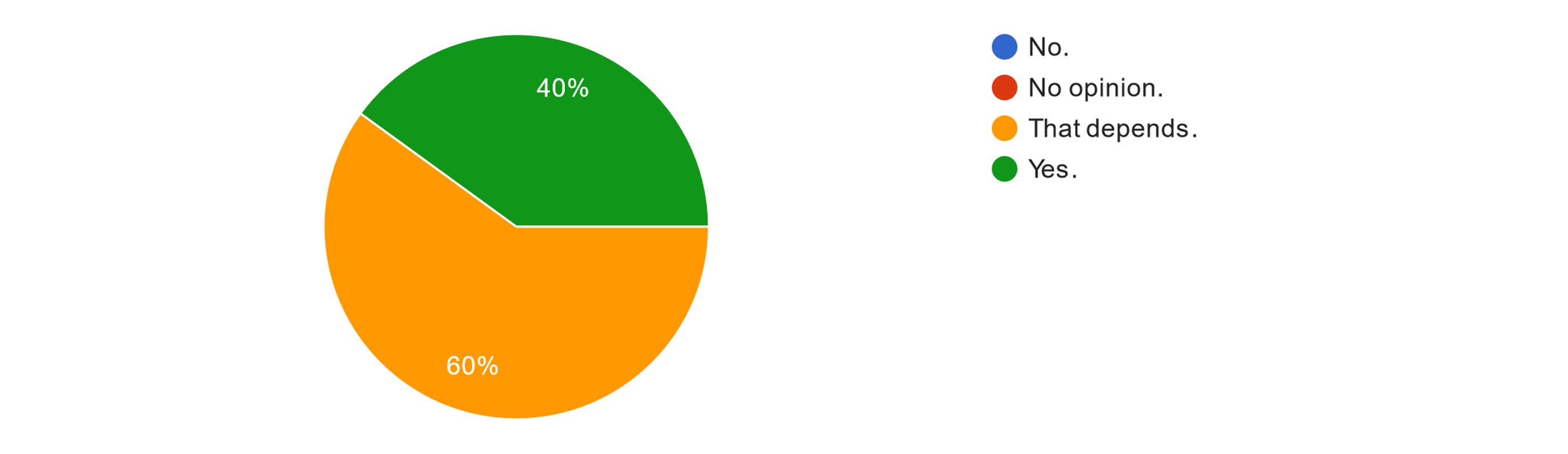}
    \caption{Question 5.}
    \label{fig:questionnaire5}
\end{figure}

   \textbf{Comments:} \begin{itemize}
    \item Some people might also find it interesting but it is not integral to the courses in my opinion.
    \item Yeah because if every student would have a module for Lean at university I think there would be more motivation to consult Lean (as a help) for a proof instead of using a website like stackexchange an coping a proof.
    \item Unfortunately, Lean is a complicated programming language (compared to Python or other imperative languages). While the basics of the tactic language can be taught in a short amount of time, when students need to push outside of those bounds or understand parts of the language itself (idea of lambda calculus, etc.), good knowledge of functional programming and a bit of type theory is required. These things are necessary when formalizing harder theorems, most of the time. So if mathematics is taught with Lean, it should be taught in conjunction with a course that teaches the language itself. Additionally, some things in mathlib might be too complicated for people just starting out with mathematics (for example that there is no extensive notion of a vector space there, everything is formulated in terms of semi-modules).
    \item I think it should be an option after the gymnasium level, and maybe it could even be a nice optional subject at the (mathematical) gymnasium. But I do not think that it should be mandatory.
    \item Yes, I think it would encourage students to work on their proofs more than currently (since you don’t need to be able to do proofs to actually pass analysis/linear algebra).
\end{itemize}

\subsubsection*{Open question: Do you plan to continue working with Lean during your studies? Why or why not?}

   \textbf{Comments:} \begin{itemize}
    \item I hope so when I have more time.
    \item Probably not often but I can see myself using it in some instances where I think that a Lean proof would help my understanding of a certain mathematical proof.
    \item Yes, I enjoy it a lot and it generally benefits my understanding of ideas and concepts. It's a lot of fun and it creates even more 'click' moments than just pen and paper (and hence a dopamine rush :D )
    \item No I do not think so, as I prefer not working with the computer and I want to learn and be able to do the proofs by 'hand'.
    \item Yes, I will try to. For example for things based on set theory this seems very effective. It’s easily understood and you get feedback whenever you’re stuck.
\end{itemize}

\subsubsection*{Open question: Would you have done something differently in the way I held my meetings with you?}

   \textbf{Comments:} \begin{itemize}
    \item For exercises that needed some specific theorem to solve it would sometimes have been useful to have a hint like 'use [theorem\_something] in the proof' as a comment because not knowing those was often the part that got me stuck and had me just looking at the solutions.
    \item No I think the meetings where always something to look forward also regarding the very relaxed and open atmosphere during them. The structure was also very good and always matched with what was happening in ``Foundations of Mathematics''.
    \item It's hard to really find fundamental improvements for the meetings with the little amount of time we had. I'm not sure how it was with the other students, but for me it was a nice experience, so I do not really have anything to improve.
    \item Maybe a touch more theory, would have helped me for example.
    \item I don’t think so, but looking back I think if the lectures started off with constructing the natural numbers (such that we could’ve started with the NNG) it might’ve been easier to get into it.
\end{itemize}

\section{Discussion}\label{sec:discussion}

   We discuss some general results from the previous section here. This chapter also contains the first author's experiences with both learning and teaching Lean. For each, we also provide some pros and cons. For the latter part, we mostly use the first-person singular perspective and write in the past tense. 

\subsection{Interviews} \label{subsec:discussionstudents}

In this section, we discuss general observations from the interviews that are not directly relevant to the results. It notes that Non-Lean students seemed more nervous than Lean students, possibly because the Lean students were familiar with the interviewer from prior meetings. While Lean might contribute to confidence, this conclusion is approached cautiously.

Not all criteria, such as the use of mathematical symbols, logic, or examples, were consistently better performed by Lean students, suggesting they benefit from using Lean but need more practice to refine these skills. The use of examples was generally weak, as it often is in proofs, though encouraging students to use examples in difficult proofs could be helpful.

Additionally, students who stopped attending meetings reported they lacked time to learn Lean due to its perceived low relevance to their studies. The importance of integrating theorem provers into the curriculum is emphasized, though no criticism is directed at students who could not attend, as the first semester is known to be stressful. It is recommended that students consider learning theorem provers later in their studies.


\subsection{Pros and Cons of Learning Lean} \label{subsec:proconlearning}

Here we briefly explain why Lean is chosen over other interactive theorem provers such as Agda, Coq, or Isabelle. Agda lacks proof automation, and Isabelle does not support dependent types, making them less suitable for some proofs. Lean and Coq are similar, but Lean has gained popularity among mathematicians due to its user-friendly features and active community.

Lean is beneficial for structuring mathematical proofs and encourages thoroughness in covering all cases, as demonstrated by interviews. However, some challenges include the temptation to rely on automated tactics without fully understanding each step and the difficulty in finding the right theorems, especially for beginners.

While Lean 3 had ample teaching materials, Lean 4 initially lacked resources during its development. However, teaching materials have since emerged, making Lean more accessible to new learners, with a dedicated webpage for contributions from the community.

\subsection{Studies with Students - Teaching (with) Lean}\label{subsec:studygroups}

For this work, students were taught how to prove mathematical statements in Lean. We chose participants from the ``Foundations of Mathematics'' course at UZH. Of the 78 students, 12 showed interest, but only seven attended, with two dropping out for various reasons. The remaining five students, all first-semester mathematics students, continued with the sessions.

Each session, lasting 45 minutes, introduced tactics and theorems needed to solve the course's exercise sheets and some extra questions, following the teaching methods outlined in Section \ref{subsec:teachingmethods}. Furthermore, one student had prior knowledge of Lean, so additional meetings were held for her to focus on more challenging proofs. In general, more than 45 minutes per week could have allowed for more explorative teaching, balancing instruction and student-led learning. Explorative teaching, however, requires significant time and effort, which students may not have been able to manage due to their other university commitments.

Initially, a document was written to guide students in installing Lean 3, but support for Lean 3 had ceased by the time sessions began \cite{ref_Lean3nothosted}. Students then switched to Lean 4, causing some delays. After adjusting to Lean 4, the sessions continued smoothly. Other than the initial software installation issue, no significant problems arose, and media-based teaching was positive.

A total of eleven meetings were planned, but only ten were held. To address student overwhelm from new material in week 11, a repetition level was prepared. The students studied natural numbers using the Natural Number Game \cite{ref_NNG, ref_NNG3}. The repetition session was postponed to the final meeting of the semester.

The attendance rate was about 80\%, with students providing reasons when they could not attend. The atmosphere was always positive, and the students' voluntary participation was appreciated.

In the next step, the study group’s performance was compared to the rest of the class. Permission was obtained to compare the exercise sheet achievements and exam marks of the entire class. Interviews with both participants and non-participants were conducted, with the results discussed in Sections \ref{subsec:interviews} and \ref{subsec:interviewresults}.

\subsection{Pros and Cons of Teaching (with) Lean} \label{subsec:proconteaching}

Teaching mathematics with Lean offers notable benefits, as observed in the sessions and interviews discussed in Section \ref{subsec:interviewresults} and supported by other sources \cite{thoma2021learning}. Lean encourages structured and precise proof construction while providing an interactive, engaging experience for students. This results in better performance in mathematical proofs.

However, Lean is not currently suitable for high school or undergraduate university teaching. High school mathematics curricula in Switzerland are already crowded, and introducing abstract topics like first-order logic or natural induction would push the limits. At the university level, teaching methods are heavily lecture-based, making it difficult to incorporate Lean without substantial restructuring. Only a few professors, like Kevin Buzzard and Patrick Massot, have integrated Lean into their teaching \cite{ref_Leanteachers}.

The challenge is not whether Lean is good for teaching but whether it aligns with the current curriculum. During the sessions, students were often overwhelmed by the number of theorems required to complete exercises. For example, in manipulating equations, several lines of code were needed, as seen in the following:

\begin{lstlisting}[backgroundcolor = \color{light-gray}, escapeinside=``]
/-Use the method of direct proof to prove the following statements.
Let x, y ∈ ℝ. If x^2 + 5y = y^2 + 5x, then x = y or x + y = 5.-/
example (x y : ℝ) :  (x ^ 2 + 5 * y = y ^ 2 + 5 * x) → 
((x = y) ∨ (x + y = 5)) := by
intro h
rw [← sub_eq_zero] at h
rw [← sub_eq_zero] at h
rw [← sub_sub] at h
rw [add_comm] at h 
rw [← add_sub] at h
rw [add_comm] at h
rw [_root_.sq_sub_sq] at h
rw [sub_zero] at h
rw [sub_eq_add_neg] at h
rw [add_assoc] at h
rw [add_comm (5*y) (-(5*x))] at h
rw [← add_assoc] at h
rw [← sub_eq_add_neg] at h
rw [sub_add] at h
rw [← mul_sub] at h
rw [← sub_mul] at h
rw [ _root_.mul_eq_zero] at h
cases' h with h1 h2
right
rw [sub_eq_zero] at h1
exact h1
left
rw [sub_eq_zero] at h2
exact h2
done
\end{lstlisting}~

While this could be solved using the \textit{have} tactic, tactics like \textit{ring\_nf} are still required. Furthermore, students faced difficulties with certain exercises, such as proving the \textit{Cantor-Bernstein-Schröder} Theorem. On paper, students could finish the proof by citing the theorem, but in Lean, they had to type a specific command:

 \begin{lstlisting}[backgroundcolor = \color{light-gray}]
exact Function.Embedding.schroeder_bernstein f_5_3_is_injective f_5_3_inv_is_injective
done
\end{lstlisting}

This made the proof less intuitive, especially if students did not know the exact command for the theorem. Therefore, students must learn both mathematics and Lean itself. To avoid overwhelming students, a hybrid approach might be more suitable. Some parts of proofs could be completed in Lean, while others could be done on paper. For example, in the Cantor-Bernstein-Schröder theorem, students could prove injectivity in Lean and finish the proof on paper. This approach could combine the best aspects of both methods, facilitated by the \textit{sorry} command in Lean.

\begin{center}
    \subsection*{Acknowledgments}
\end{center}
We acknowledge the invaluable, ongoing support of the Lean community. Special thanks to Marius Furter for his clear explanations of the logic behind Lean tactics in \cite{mattiamaster}, and to the Institute of Mathematics at the University of Zurich for enabling this study. Lastly, sincere thanks to all the students involved in this research; their dedication and collaboration were essential to making this study possible.

\medskip

\bibliography{References}

\end{document}